\documentstyle[amssymb]{article}

\newcommand{\e}{\epsilon}
\newcommand{\be}{\bar{\epsilon}}

\newcommand{\Lap}{\Delta}
\newcommand{\Lbe}{L_{\bar{\e}}}
\newcommand{\Lbel}{L_{\bar{\e}^{(\ell)}}}
\newcommand{\Gbe}{G_{\bar{\e}}}

\newcommand{\bue}{\bar{u}_{\bar{\epsilon}}}
\newcommand{\fbe}{f_{\bar{\epsilon}}}
\newcommand{\del}{\partial}
\newcommand{\OS}{\Omega \setminus \Sigma}
\newcommand{\el}{\e^{(\ell)}}
\newcommand{\bel}{\be^{(\ell)}}
\newcommand{\ubel}{u_{\bel}}
\newcommand{\vbel}{v_{\bel}}
\newcommand{\fbel}{f_{\bel}}
\newcommand{\tubel}{\tilde u_{\be^{(\ell)}}}
\newcommand{\tvbel}{\tilde v_{\be^{(\ell)}}}
\newcommand{\tfbel}{\tilde f_{\be^{(\ell)}}}
\newcommand{\bubel}{\bar{u}_{\be^{(\ell)}}}
\newcommand{\bvbel}{\bar{v}_{\be^{(\ell)}}}
\newcommand{\bfbel}{\bar{f}_{\be^{(\ell)}}}
\newcommand{\RR}{\Bbb R}
\newcommand{\RNO}{\Bbb R^N \setminus \{0\}}
\newcommand{\Rnk}{\Bbb R^n \setminus \Bbb R^k}
\newcommand{\Ts}{\cal T_{\sigma}}
\newcommand{\C}{\cal C}
\newcommand{\Obe}{\Omega_{\be}}
\newcommand{\Obel}{\Omega_{\be^{(\ell)}}}
\newcommand{\Sig}{\Sigma}
\newcommand{\LL}{\Bbb L}

\begin{document}

\title{A construction of singular solutions for a semilinear \\
elliptic equation using asymptotic analysis}
\author{Rafe Mazzeo \\ Stanford University \thanks{Supported by
the NSF under Grant DMS-930326 and a Young Investigator Fellowship,
and by a Sloan Foundation Fellowship}
\and Frank Pacard \\ ENPC-Cergrene}

\maketitle

\newtheorem{theorem}{Theorem}
\newtheorem{lemma}{Lemma}
\newtheorem{proposition}{Proposition}
\newtheorem{definition}{Definition}
\newtheorem{corollary}{Corollary}
\newtheorem{remark}{Remark}

\begin{abstract}
The aim of this paper is to prove the existence of weak solutions
to the equation $\Delta u + u^p = 0$ which are positive in
a domain $\Omega \subset {\Bbb R}^N$, vanish at the boundary, and have
prescribed isolated singularities. The exponent $p$ is required to
lie in the interval $(N/(N-2), (N+2)/(N-2))$. We also
prove the existence of solutions
to the equation $\Delta u + u^p = 0$ which are positive in
a domain $\Omega \subset {\Bbb R}^n$ and which are
singular along arbitrary smooth $k$-dimensional submanifolds in the
interior of these domains provided $p$ lie in the interval
$((n - k)/(n-k-2), (n-k+2)/(n-k-2))$. A particular case
is when $p = (n+2)/(n-2)$, in which case solutions correspond
to solutions of the singular Yamabe problem. The method used
here is a mixture of different ingredients used by both
authors in their separate constructions of solutions to the
singular Yamabe problem, along with a new set of scaling techniques.
\end{abstract}

\tableofcontents

\section{Introduction and statements of main results}

In this paper we construct solutions with prescribed singularities
for the semilinear elliptic equation $\Delta u + u^p = 0$
and other closely related equations, for a certain range of values
of the exponent $p$, in a variety of situations. The solutions will
have singularities prescribed along a disjoint union of submanifolds
of varying dimension. We now describe our results, starting with the
simplest case, when the prescribed singular set is discrete.

Suppose $\Omega$ is any bounded open set in $\Bbb R^N$, with
smooth boundary. Consider the equation
\begin{equation}
\left\{ \begin{array}{rllll}
          -\Delta u  & = & u^{p} & \mbox{in} & \Omega \\
                  u  & = & 0           & \mbox{on} & \del  \Omega ,
        \end{array}
\right.
\label{eq:delta-u=u-p}
\end {equation}
A weak solution of (\ref{eq:delta-u=u-p}) is a function $u$ which solves
this equation on all of $\Omega$ in the sense of
distributions. In particular, any such solution must belong to
$L^{p}(\Omega)$.  The singular set of a weak solution $u$,
$\mbox{sing\,}(u)$, is the complement in $\Omega$ of the set of
points where $u$ is continuous, and hence smooth.

Our first result concerns the case when the singular set is finite.
We single it out since its proof is slightly simpler than the more general
case.
\begin{theorem}
Let $\Omega$ be as above, and suppose that $\Sigma = \{x_1, \dots,
x_K\} \subset \Omega$ is any finite set of points. Suppose also that
the exponent $p$ lies in the range
\begin{equation}
{{N \over {N-2}}} \leq  p < {{N+2}\over {N-2}}.
\label{eq:range-of-p-is}
\end{equation}
Then there is a $K$ parameter family of positive weak solutions $u$ of
(\ref{eq:delta-u=u-p}) with $\mbox{sing\,}(u) = \Sigma$. In fact, the
solution space of this equation is locally a $K$ dimensional real
analytic variety.
\label{th:ex-is-sing}
\end{theorem}

It is known that if the exponent $p$ is less than $N/(N-2)$, then
any weak solution of (\ref{eq:delta-u=u-p}) must be smooth on
all of $\Omega$. The existence of solutions of this equation
with prescribed isolated singularities when $p$ lies in
the interval $N/(N-2) \le p < p_0$, where $p_0$ is some value
close to $N/(N-2)$ (and in particular, less than $(N+2)/(N-2)$),
has already been solved by the second author in \cite{P1} and in \cite{P2}.
When $p = (N+2)/(N-2)$,
which is the so-called critical exponent, the problem becomes
conformally invariant. There is then a loss of compactness and the problem
consequently becomes much more difficult.
Solutions now correspond to metrics of constant positive scalar curvature
which are complete in a neighbourhood of the singular points. It is
geometrically more natural in this case to replace the domain by $S^N$
(or in fact, any other compact manifold of positive scalar
curvature); the operator $\Delta$ then needs to be replaced
by the conformal Laplacian $\Delta - \frac{N-2}{4(N-1)}R_0$,
where $R_0$ is the scalar curvature of the background manifold.
An additional source of difficulties in this case is
that the position of the singularities is no longer
necessarily arbitrary. The general existence result for this geometric
problem (when the background manifold is the sphere) was obtained by
R.\  Schoen \cite{S}. The recent work of the first author, along with
D.\  Pollack and K.\  Uhlenbeck \cite{M-P-U} examines the moduli space of
solutions for this problem.

In any event, Theorem~\ref{th:ex-is-sing} extends the result
of \cite{P1} and \cite{P2} to the full range of subcritical exponents.

It is also possible to prove existence of solutions of (\ref{eq:delta-u=u-p})
which are singular along submanifolds of higher dimension.
Now let $\Sigma = \cup_{i=1}^K \Sigma_i$, where each $\Sigma_i
\subset \Omega$ is a smooth submanifold without boundary of dimension
$k_i \ge 0$. We shall suppose that $\Omega \subset \RR^n$, since we
reserve the symbol $N$ (or more properly $N_i = n-k_i$) for the dimension
of the normal space to each $\Sigma_i$.
\begin{theorem}
Suppose that the exponent $p$ satisfies
\[
\frac{n-k_i}{n-k_i - 2} \leq  p < \frac{n-k_i + 2}{n-k_i - 2},
\label{eq:range-of-p-gen}
\]
or in other words
\[
n - \frac{2p+2}{p-1} < k_i \leq n - \frac {2p}{p-1},
\]
for $i = 1, \dots, K$. Then there is a positive weak solution $u$ of
(\ref{eq:delta-u=u-p}) with $\mbox{sing\,}(u) = \Sigma$.  Provided
at least one of the $k_i > 0$, there is an infinite dimensional
space of solutions for this problem.
\label{th:ex-gen-sing}
\end{theorem}

Again, the existence of solutions for this equation with higher
dimensional singular sets of dimension $k$ has been solved by
Y.\ Reba\"{\i} \cite{YR} when the exponent $p$ lies in the
interval $(n-k)/(n-k-2) \le p < p_k$, where $p_k$ is close to
the lower endpoint of this interval, and in particular, less
than $(n-k + 2)/(n-k-2)$.  We shall only prove this result when
$p > (n-k_i)/(n - k_i - 2)$ since it would complicate the notation
to cover the borderline case where $p$ attains the lower limit of
this interval of values. However, the arguments below could be extended
without undue difficulty to cover this value as well.

A special case of Theorem~\ref{th:ex-gen-sing} is the singular
Yamabe problem. This occurs when the exponent $p$ attains
the critical value $(n+2)/(n-2)$. In this case, the dimensions
of the $\Sigma_i$ are allowed to lie in the range $0 < k_i < (n-2)/2$.
As discussed above, it is more natural now to replace $\Omega$ by the
round sphere $(S^n,g_0)$, or indeed by any compact manifold $(M,g_0)$ with
(constant) nonnegative scalar curvature, and to replace the Laplacian
by the conformal Laplacian $L_0 = \Delta - \frac{n-2}{4(n-1)}R(g_0)$.
The relevant equation now is
\begin{equation}
-L_0 u =  u^{\frac{n+2}{n-2}}
\label{eq:sing-Yam}
\end{equation}
on $M$. At least when the background manifold is the sphere, it follows
from work of R.\ Schoen and S.\ T.\ Yau \cite{SY} that solutions to this
singular Yamabe problem (for which the corresponding conformally
related metrics $g = u^{4/(n-2)}g_0$ are complete and have bounded
Ricci curvature) exist only if the dimension of the singular set
is less than or equal to $(n-2)/2$. Thus, Theorem~\ref{th:ex-gen-sing}
implies
\begin{theorem}
Let $(M,g_0)$ be any compact manifold with constant nonnegative scalar
curvature. Let $\Sigma \subset M$ be any finite disjoint union
of smooth submanifolds $\Sigma_i$ of dimensions $k_i$ with
$0 < k_i \leq \frac{n-2}2$.  Then there is an infinite dimensional
family of complete metrics on $M \setminus \Sigma$ with constant
positive scalar curvature.
\label{th:sing-Yam}
\end{theorem}
Together with the result of \cite{P}, which treats the case
$k_i=\frac{n-2}{2}$ when $n$ is even, this theorem settles the
question of existence of solutions to the singular Yamabe problem
(with constant positive scalar curvature) whenever $\Sigma$ is a
finite disjoint union of smooth submanifolds with dimensions greater
than zero, but less than or equal to $(n-2)/2$.

Finally, we may apply the arguments of \cite{M-S} to slightly
refine Theorem \ref{th:ex-gen-sing}:
\begin{theorem} Let $\Sigma$ be any finite disjoint union of
$\C^{3,\alpha}$ submanifolds in $\Omega$ of dimensions $k_i$ satisfying
the restrictions above. Then (\ref{eq:delta-u=u-p}) has an infinite
dimensional family of solutions.
\end{theorem}

The basic idea of Theorems~\ref{th:ex-gen-sing} and~\ref{th:sing-Yam}
is that solutions of (\ref{eq:delta-u=u-p}) with positive
dimensional singular sets may be obtained as perturbations
of solutions to the problem on the fibres of the normal
bundle of $\Sigma$. The exponent $\frac{n+2}{n-2}$ is subcritical
for the induced problem on these fibres, which are of lower dimension
than the ambient space. In particular, the problem with critical exponent
may be reduced to a subcritical one.

The paper is organized as follows. First we analyze the
asymptotics of rotationally invariant solutions of the
problem on $\Bbb R^N \setminus \{0\}$.  These are then
used to construct approximate solutions for each of the problems.
At this point we then give a more detailed outline of the strategy
of the proof. The main task is to analyze the
linearization of the nonlinear operator around these approximate
solutions, and to prove that these linearizations are surjective
on appropriate function spaces. This occupies the bulk of the
paper. After this is accomplished, the exact solutions are
obtained by a fixed point argument. The deformation spaces for the
solutions of these operators, and the arguments necessary to
replace $\Omega$ by a manifold $M$ are discussed in the last section.

The authors are grateful to Dan Pollack for having carefully read
the manuscript and making a number of suggestions to improve
the exposition.

\section{Singular radial solutions on $\Bbb R^N \setminus \{0\}$}

In this section, we recall some well known facts which appear for
example in \cite{C-G-S}.

\begin{proposition}
For any exponent $p \in ({N \over {N-2}}, {{N+2}\over {N-2}})$, there exists a
one parameter family  of weak solutions $ u_{\e},\ \e >0$, for the equation
\begin{equation}
-\Delta u = u^p \quad \mbox{in} \quad {\Bbb R}^{N},
\label{eq:nonlinear-p}
\end{equation}
such that the $u_{\e}$ are radial, singular only at the origin, and
satisfy the following properties~:
\begin{itemize}

\item $u_{\e}(r) > 0$ for $0 < r < \infty$.

\item $u_{\e}$ can be written
as
\[
u_{\e}(x)\equiv |x|^{-{2\over {p-1}}}v_{1}(-\log(|x|/\e)),
\]
where the function $v_{1}$ is bounded independently of $\e$.

\item $\lim_{t\rightarrow +\infty} v_{1}(t) = v_{\infty} >0,$ where
\begin{equation}
v_{\infty}^{p-1} = {2 \over {p-1}} \left( N - {{2p}\over {p-1}}\right),
\label{eq:v-infty}
\end{equation}
and in particular, is independent of $\e$.

\item  $\lim_{t\rightarrow -\infty} e^{-t(N-\frac{2p}{p-1})} v_{1}(t)
<+\infty$; in particular, for $|x|$ large
\[
u_\e (x) =c (\e^{N -{{2p}\over {p-1}}} +O(|x|^{-1}) |x|^{2-N},
\]
 where the constant $c$ is independent of $\e$.

\item Finally, $||v_{1}||_{L^{\infty}}^{p-1} < {{p+1}\over 2}
v_{\infty}^{p-1}.$
\end{itemize}
\label{pr:base}
\end{proposition}
{\bf Proof~:}
First define a new independent variable $t = -\log |x|$ and set
\begin{equation}
u(x)\equiv |x|^{-{2\over {p-1}}}v(-\log(|x|)).
\label{eq:change-u-v}
\end{equation}
Then the new function $v(t)$ satisfies
\begin{equation}
\del_{tt} v -\left( N  - {{p+1}\over {p-1}}\right)\del_t v
-{{2p} \over {p-1}} \left( N  - {{2p}\over {p-1}}\right)v +v^{p}=0.
\label{eq:radial-p}
\end{equation}
Now look at the phase-plane portrait for this equation in
the $(v, v_{t})$ plane. The two equilibrium points
are $(0,0)$ and $(v_{\infty},0)$, where $v_{\infty}$ is defined in
(\ref{eq:v-infty}); the first of these is a saddle point
and the second a stable equilibrium. There is a single orbit
issuing from $(0,0)$ and tending to $(v_{\infty},0)$
as $t \rightarrow \infty$.  Let $v_1(t)$ be one of the functions corresponding
to this orbit; it is determined only up to the choice of its initial
Cauchy data $(v_1(0),v_1'(0))$.
Now let $v_\e(t) = v_1(t + \log \e)$. As $\e$ varies over $(0,\infty)$
we obtain all solutions of (\ref{eq:radial-p}) which converge to
$0$ as $t \rightarrow -\infty$ and to $v_{\infty}$ as $t \rightarrow
\infty$.

We now obtain an upper bound on the function $v_1(t)$.
It is easy to see that the trajectory $(v_1(t),v_1'(t))$
is contained within the homoclinic orbit of the Hamiltonian system
\begin{equation}
\del_{tt} w
-{{2p} \over {p-1}} \left( N  - {{2p}\over {p-1}}\right)w +w^{p}=0
\label{eq:hamilton-p}
\end{equation}
which tends to $(0,0)$ as $t$ tends both to $+\infty$ and $-\infty$.
Let $(w_1(t),w_1'(t))$ parametrize this orbit. Then we conclude that
\[
\sup\, v_1 \leq \sup\, w_1.
\]
The conservation of Hamiltonian energy for (\ref{eq:hamilton-p})
now shows that
\begin{equation}
{1 \over 2} (\del_{t} w_1 )^2
-{{p} \over {p-1}} \left( N  - {{2p}\over {p-1}}\right)w_1^2 +{1 \over
{p+1}}w_1^{p+1}=0.
\end{equation}
$w_1$ attains its supremum when $w_1' = 0$, so
we obtain the upper bound
\[
v_1^{p-1} (t) < {{p+1} \over 2} v_{\infty}^{p-1} \qquad \mbox{for all}
\qquad t \in {\Bbb R}.
\]
This ends the proof of the Proposition.\hfill $\Box$
\begin{remark}
For any $\e >0$, whenever $u(x)$ is a solution of
(\ref{eq:nonlinear-p}), $\e^{-{{2}\over {p-1}}}u(\e^{-1}\,x)$
is another solution of this same equation. This dilation invariance
corresponds exactly to the translation invariance of the equation
(\ref{eq:radial-p}) for $v$. We let $u_1(x)$ be the solution
corresponding to $v_1$ and
\[
u_{\e}(x)\equiv \e^{-{{2}\over {p-1}}}u_{1}( \e^{-1}\,x)
\]
\end{remark}
Later we will also need the following~:
\begin{remark}
By this rescaling we can always assume that, for any given constant
$\alpha >0$, $v_{1}$ may be chosen so that
\begin{equation}
\sup_{t \geq 0} v_{1}(t) \leq \alpha.
\end{equation}
\label{re:small-outside}
\end{remark}
The final remark relates the construction of approximate solutions
here to those used in \cite{M-S}:
\begin{remark}
The stable point $(v_{\infty},0)$ corresponds to the singular solution
\[
u_0 (x) \equiv v_{\infty} |x|^{-{2\over {p-1}}}.
\]
This solution is invariant by the scaling described above. The
fact that $u_{1}$ is not dilation invariant is important
in the construction of approximate solutions.
\end{remark}

\section{Function spaces}

In this section we define the weighted H\"older spaces
$\C^{k,\alpha}_{\nu}(\OS)$ appropriate for
this problem. Roughly speaking, the functions in these spaces
are products of powers of the distance to $\Sigma$ with functions
whose H\"older norms are invariant with respect to scaling by
dilations from any point on $\Sigma$.

We shall use local Fermi coordinates around each component of $\Sig$
to define these spaces. When $\Sig_i$ is a point,
these are simply (geodesic) polar coordinates $(r,\theta)$ around that
point. When $\Sig_i$ is higher dimensional, let $\Ts^{(i)}$ be the tubular
neighbourhood of radius $\sigma$ around $\Sigma_i$. It is well known
that $\Ts^{(i)}$ is a disk bundle over $\Sigma_i$, and is diffeomorphic
to the disk bundle of radius $\sigma$ in the normal bundle
$N\Sigma_i$. Using the metric, this diffeomorphism is canonical.
The Fermi coordinates in this tubular neighbourhood will be constructed
as coordinates in the normal bundle, transferred to $\Ts^{(i)}$ via this
fixed diffeomorphism. Here $r$ is the distance to $\Sigma_i$, which is well
defined in $\Ts^{(i)}$ and smooth away from $\Sigma_i$
provided $\sigma$ is small enough, $y$ is a local coordinate system on
$\Sigma_i$, and $\theta$ is the angular variable on the sphere in
each normal space $N_y\Sigma_i$.  Let $B_{N,\sigma}$ denote
the ball of radius $\sigma$ in $N_y\Sigma_i$. We shall let $x$ denote
the rectangular coordinate in these normal spaces, so that
$r = |x|$ and $\theta = x/|x|$.

Let us also fix a function $\rho > 0$ in $\C^{\infty}(\OS)$ with $\rho$
equal to the polar distance $r$ in each $\Ts^{(i)}$.
Let $w$ be a function in this tubular neighbourhood, and define
\[
||w||_{0,\alpha,0}^{\Ts^{(i)}}=
\sup_{z \in \Ts^{(i)}} |w| + \sup_{z,\tilde{z} \in \Ts^{(i)}}
\frac{ (r + \tilde{r})^{\alpha}|w(z) - w(\tilde{z})|}{
|r - \tilde{r}|^{\alpha} + |y - \tilde{y}|^{\alpha} +
(r+\tilde{r})^{\alpha} |\theta - \tilde{\theta}|^{\alpha}}.
\]
Here $z, \tilde{z}$ are any two points in $\Ts^{(i)}$ and
$(r,y,\theta),\ (\tilde{r},\tilde{y}, \tilde{\theta})$ are
their Fermi coordinates.

\begin{definition}
The space $\cal C^{k,\alpha}_{0}(\OS)$ is defined to be the set of all
$w \in \C^{k,\alpha}(\OS)$ for which the norm
\[
||w||_{k,\alpha,0} \equiv ||w||_{k,\alpha, \Omega_{\sigma/2}}
+ \sum_{i=1}^K \sum_{j = 0}^k || \nabla^j w||_{0,\alpha}^{\Ts^{(i)}}
\]
is finite. Here $\Omega_{\sigma/2} = \Omega \setminus \cup_{j=1}^K
\cal{T}_{\sigma/2}^{(i)}$. Then, for any $\gamma \in \RR$,
\[
\C^{k,\alpha}_{\gamma}(\OS) = \{w = \rho^{\gamma} \bar{w}: \bar{w} \in
\C^{k,\alpha}_{0}(\OS)\}.
\]
\end{definition}
Thus functions in $\cal C^{k,\alpha}_{\gamma}(\OS)$
are allowed to blow up like $\rho^{\gamma}$ near the $\Sig_i$. Note that
when $\Sig_i$ is positive dimensional, functions in $\C^{k,\alpha}_{\gamma}$
may be differentiated in the `tangential' direction only at the expense
of giving up a power of $\rho$. Equivalently, their derivatives
with respect to up to $k$-fold products of the vector fields
$r\del_r$, $r\del_y$, $\del_{\theta}$ blow up no faster than $\rho^{\gamma}$.

We collect a few essentially trivial remarks about these spaces~:
\begin{lemma}
\begin{enumerate}
\item If $\gamma \in (-N,0)$ then ${\cal C}_{\gamma}^{k, \alpha }{(\OS)}
\subset L^{p}(\Omega )$ for $p \in (1, -N/\gamma)$.
\item For $w \in {\cal C}_{\gamma}^{k, \alpha}$ and $v\in {\cal
C}_{\gamma'}^{k, \alpha}$, then
$wv \in {\cal C}_{\gamma +\gamma'}^{k, \alpha}$ and also
\[
 ||wv||_{k,\alpha,\gamma+\gamma'} \leq 4
||w||_{k,\alpha,\gamma} ||v||_{k,\alpha,\gamma'}.
\]
\item Given $w\in {\cal C}_{\gamma}^{k, \alpha}$ with $w \geq 0$
and if $p > k+1$, then $w^{p}\in {\cal C}_{p\gamma }^{k, \alpha}$ and
\[
 ||w^{p}||_{k,\alpha,p\gamma} \leq c_p ||w||_{k, \alpha ,\gamma}^{p},
\]
for some constant $c_{p}>0$ only depending on $p$.

\item If $ w\in {\cal C}^{k+1,0}_\gamma$ and
$ |\nabla w| \in {\cal C}^{k,0}_{\gamma-1}$, then for any $\alpha \in
(0,1)$, $w\in {\cal C}_{\gamma}^{k, \alpha}$.
\end{enumerate}
\label{le:fcnspacefacts}
\end{lemma}

We shall also need a slightly better estimate for the norm of $w^p$.
\begin{lemma}
For any $\eta>0$, $\gamma \in \RR$, there exists some $\theta>0$, depending
only on $\gamma$ and $\eta$, such that if $||w||_{0,\alpha,\gamma} <
\theta$ then
\[
||w^p||_{0, \alpha,\gamma -2 } \leq \eta ||w||_{ 0, \alpha ,\gamma }.
\]
\label{le:v^p}
\end{lemma}
{\bf Proof~:} The proof relies on the simple fact that
\[
||w^p||_{ 0, \alpha , p\gamma } \leq c_{}
\eta^{p-1} ||w||_{ 0, \alpha ,\gamma}.
\]
Since $ \gamma > -{2 \over {p-1}}$, we see that $p\gamma \geq  \gamma -2$, and
the result follows immediately.
\hfill $\Box$

\begin{definition}
\[
\C^{k,\alpha}_{\gamma,\cal D}(\OS) = \{w \in \C^{k,\alpha}_{\gamma}(\OS):
w = 0 \mbox{ on } \del \Omega\}.
\]
\end{definition}

Finally, we shall also need to use H\"older spaces on $\RNO$,
respectively $\Rnk$, which have different decay properties at the origin,
respectively at $\RR^k$, than at $\infty$.  Now we simply use the global
polar coordinates $(r,\theta)$ in $\RNO$, and cylindrical
coordinates $(r,\theta,y)$ in $\Rnk = \RNO \times \RR^k$.
\begin{definition}
For any $\gamma, \gamma' \in \Bbb R$, the space $\cal C^{k,\alpha}_{\gamma,
\gamma'}(\RNO)$ consists of all functions $w$ for which the norm
\[
|| w ||_{k,\alpha,\gamma,\gamma'} \equiv \sup_{ B_2(0)}
||w||_{k,\alpha,\gamma} +
\sup_{\Bbb R^N \setminus B_1(0)} ||w||_{k,\alpha,\gamma'}
\]
is finite. The definition of $\C^{k,\alpha}_{\gamma,\gamma'}(\Rnk)$
is similar; we need only replace $B_j(0)$ by the tube of radius
$j$ around $\RR^k$, $j =1,2$.
\end{definition}

Although these H\"older spaces are our primary tools, we shall also
need to refer on occasion to a family of weighted $L^2$ spaces
and their associated Sobolev spaces. These will be needed
only for nonisolated singularities, so we restrict ourselves
to that case.
\begin{definition}
The weighted space $r^{\delta}L^2$ on $\Rnk$ is defined by
\[
r^{\delta}L^2(\Rnk; r^{N-1}\,dr\,d\theta\,dy) =
 \{ w \in L^2_{\text{loc}}(\Rnk): \int |w|^2 \,
r^{N-1-2\delta}\,dr\,d\theta\,dy < \infty\}.
\]
The space $\rho^{\delta}L^2(\OS)$ (relative to standard Euclidean volume
measure) is defined similarly, using the function $\rho$ above.
\end{definition}
There are associated weighted Sobolev spaces $r^{\delta} H^s_e(\Rnk)$
when $s$ is a positive integer, where the subscript $e$ signifies
that these spaces are defined with respect to differentiations by
the vector fields $r\del_r$, $r\del_y$ and $\del_\theta$. (When
$s$ is an arbitrary real number, they may be defined by duality and
interpolation.) Note that
\[
\C^{k,\alpha}_{\gamma,\sigma}(\Rnk) \subset r^{\delta}L^2(\Rnk)
\quad \mbox{provided}\quad \sigma + (N-2)/2 < \delta < \gamma + (N-2)/2
\]
\begin{equation}
\C^{k,\alpha}_{\gamma}(\OS) \subset \rho^{\delta}L^2(\OS) \quad
\mbox{provided}\quad \delta < \gamma + (N-2)/2.
\label{eq:space-incl}
\end{equation}

\section{Construction of approximate solutions}

Now that the family of radial solutions $u_{\e}(x)$ to
(\ref{eq:nonlinear-p}) on $\RNO$ has been
introduced, we can construct the approximate solutions for
these problems. We also derive estimates of how far these solutions
differ from exact solutions in terms of the weighted H\"older
norms of the last section.

\subsection{Approximate solutions with isolated singularities}

Approximate solutions for (\ref{eq:delta-u=u-p}) which are
singular precisely at the points of $\Sigma = \{x_1, \dots, x_K\}$
are constructed by superimposing appropriately translated and dilated copies
of the singular radial solutions. The $K$ free parameters in
the exact solution correspond to the $K$ different dilation
parameters we can prescribe independently at each of the
singular points.  First, choose $\chi(x) \in \cal
C^{\infty}_0({\RR}^N)$ with $\chi = 1$ when
$|x| \leq 1$ and $\chi = 0$ when $|x| \geq 2$, set $\chi_{R}(x) \equiv
\chi (x/R)$. Also choose
$R > 0$ with $R < R_0 \equiv \inf_{i\neq j}
( \mbox{dist} (x_{i}, x_{j}))/2$).  Let $\be = (\e_1,
\cdots, \e_K\}$ be a $K$-tuple of dilation parameters. Now
define
\begin{equation}
\bue (x) =\sum_{i=1}^{K} \chi_R
(x-x_{i})u_{\epsilon_i}(x-x_{i})
=\sum_{i=1}^{K} \chi ({{x-x_{i}}\over R})
\epsilon_i^{-{2 \over {p-1}}}u_{1}({{x-x_{i}}\over \epsilon_i}).
\label{eq:approximate-N}
\end{equation}

Set $\fbe = \Delta \bue +\bue^p$. After some computation, we get
\[
\fbe = \sum_{i=1}^{K}  u_{\e_i}(x-x_i)\Delta \chi_R (x-x_i)
+2 \nabla u_{\e_i}(x-x_i)\nabla \chi_R (x-x_i) +\sum_{i=1}^{K}
(\chi^{p}_R-\chi_R ) (x-x_i) u_{\e_i}^p (x-x_i).
\]

\begin{lemma}
For any $\gamma \in \RR$, there exists a constant $c$, depending on $R$ and
$\gamma$, such that
\[
|| \fbe ||_{0, \alpha,\gamma-2 } \leq c \e_0^{N -{{2p}\over {p-1}}},
\quad \text{provided each}\quad \e_i \le \e_0 \le 1.
\]
\label{le:estimate-1}
\end{lemma}
{\bf Proof~:}
For $|x-x_{i}|\in (R, 2R)$, we get
\[
|\bue(x)| \leq  c\, \e_0^{N -{{2p}\over {p-1}}}
\]
and
\[
|\nabla \bue(x)| \leq c\,\e_0^{N -{{2p}\over {p-1}}},
\]
for some constant $c>0$ depending on $N$ and $R$. The result follows at once.
\hfill $\Box$

\subsection{Approximate solutions in the general case}

Now suppose $\Sigma = \cup_{i = 1}^K \Sigma_i$ where each
$\Sigma_i$ is a smooth submanifold in $\Omega$ of dimension
$k_i$.  To simplify the notation here, we assume that all $k_i > 0$.

In terms of the Fermi coordinates introduced in \S3,
the Euclidean Laplacian on $\Ts^{(i)}$ can be written locally
in terms of the Laplacians for $\Sigma_i$ and
$N_y\Sigma_i$:
\begin{equation}
\Delta  = \Delta_{N} + \Delta_{\Sigma_i} +
e_{1}\cdot \nabla ^{2} +
e_{2}\cdot\nabla,
\label{eq:lap-fermi}
\end{equation}
where $e_{1}$ and $e_{2}$ satisfy
\[
||e_{1}||_{0,\alpha ,1} + ||e_{2}||_{0,\alpha ,0} \leq c,
\]
for some constant  $c >0$ which does not depend on $\alpha$, nor
on $x,y$.
This is discussed in~\cite{M-S} and \cite{FMO}.

Now choose a smooth cut-off function $\chi$ on ${\RR}^N$ as before,
which only depends on $|x|$, and such that $\chi(x) = 0$ for $|x| > 2$
and $\chi(x) = 1$ for
$|x| < 1$. Also, set $\chi_{R}(x) \equiv \chi (x/R)$.
If $0 < \e_i < 1$ and $R <\sigma /2$, then define, in some
neighborhood of $y_{0}\in \Sigma_i$, the function
\[
\bar{u}_{\e_i}(x,y) \equiv \e_i^{-\frac{2}{p-1}}
u_{1}( x / \e_i) \chi_{R} (x)\equiv u_{\e_i}(x) \chi_{R}(x).
\]
Since this function only depends on the variable in the
normal space, and is independent of the angular variable $\theta$,
it is clear that it may be defined globally on all of $\Ts^{(i)}$.
Now let $\bue = \sum_{i=1}^K \bar{u}_{\e_i}$.

As before, let
\[
f_{\be} = \Lap \bue + \bue^p.
\]
Using (\ref{eq:lap-fermi}), we compute that
\[
f_{\be} = e_{1}(x,y)\cdot \nabla^{2}\bue
+ e_{2}(x,y)\cdot \nabla \bue +
\bue\Delta_{N}\chi_{R}
+ 2\nabla \chi_{R}\nabla \bue +(\chi^{p}_{R} - \chi_{R})
\bue^{p}.
\]

\begin{lemma}
There exists some $c>0$ depending on $\gamma$ but independent of $\e_0 <1$
such that if each $\e_i \le \e_0$, then
\[
||f_{\be}||_{0, \alpha, \gamma-2 } \leq c \e_0^{q},
\]
where $q=\frac{p-3}{p-1}-\gamma$.
\label{le:Estimate-1}
\end{lemma}
{\bf Proof~:} The estimate follows at once from similar estimates
already used in Lemma~\ref{le:estimate-1} and also from the estimates
given in (\ref{eq:lap-fermi}).\hfill $\Box$

This exponent $q$ is strictly positive by our assumptions
on $p$, provided $\gamma > \frac{-2}{p-1}$.

\section{Outline of the proof}
In both the cases of solutions with isolated or more general singularities,
we wish to solve the equation (\ref{eq:delta-u=u-p}) by perturbing
the approximate solutions $\bue$ to this problem we obtained in
the preceding section. That is, we wish to represent the
exact solution $u$ as a sum $u = \bue + v$, where $v$ is small
compared to $\bue$. Thus we wish to solve
\[
\Lap(\bue + v) + (\bue + v)^p = 0,
\]
or, what is the same thing,
\begin{equation}
L_{\be}v + f_{\be} + Q(v)=0.
\label{eq:nonlin}
\end{equation}
Here
\begin{equation}
Q(v) = (\bue + v)^p - \bue^p - p\bue^{p-1}v
\label{eq:quad-remain}
\end{equation}
is the remainder term, which is quadratically small if $v$ is smaller
than $\bue$ in an appropriate H\"older norm, and
\begin{equation}
\Lbe = \Lap + p\bue^{p-1}
\label{eq:linearized-op}
\end{equation}
is the linearization about this approximate solution.

We shall show that there is some $\e_0 > 0$ and some constant $0 < c < 1$
such that for any $\e \le \e_0$, if $c \e \le \e_i \le \e$ for all $i$,
there is an exact solution of (\ref{eq:delta-u=u-p}) which is a small
perturbation of $\bue$. Writing $u = \bue + w$ as before, where $w \in
\C^{2,\alpha}_{\nu}(\OS)$,
we can ensure that $u$ is singular at $\Sigma$ provided $\nu > -2/(p-1)$.
In addition, we shall be able to prove that the solution space is a
$K$-dimensional manifold.

When any of the $k_i > 0$, then
we shall show later that the solution space is infinite dimensional.

The analysis of this linearization is the fundamental issue of this paper.
It is easy to see that
\begin{equation}
\Lbe: \C^{2,\alpha}_{\gamma}(\OS) \longrightarrow \C^{0,\alpha}_{
\gamma-2}(\OS)
\label{eq:lin-map}
\end{equation}
is bounded for any $\gamma \in {\RR}$ and $\be$.
We must show that (\ref{eq:lin-map}) is surjective for some $\gamma >
-2/(p-1)$, provided that each $\e_i$ is sufficiently small.
Once this is shown, then it follows that $\Lbe$ has a bounded right inverse
\begin{equation}
\Gbe: \C^{0,\alpha}_{\nu-2}(\OS) \longrightarrow \C^{2,\alpha}_{\nu}(\OS).
\label{eq:rinv-map}
\end{equation}
We must also show that the norm of this map is bounded independently of
the $\e_i$ for $\e_i < \e_0$.  However, for the range of values of
$\nu$ for which (\ref{eq:lin-map}) is surjective, it is not injective,
so that $\Gbe$ is not unique. As usual, a good choice for a right inverse
is the one which maps into a fixed complement of the nullspace of $\Lbe$,
for example, the orthogonal complement with respect to some weighted
$L^2$ structure.

To prove surjectivity of (\ref{eq:lin-map}) we prove both that
its range is closed and that the cokernel is trivial.  When $\Sigma$ is
a discrete set, the closedness of the range is rather elementary,
but when some component of $\Sigma$ has positive dimension this fact
is somewhat deeper. Fortunately, the construction of pseudodifferential
right parametrices for the general class of elliptic `edge operators,'
of which $\Lbe$ is a particular example, falls within the scope of
the theory developed by the first author in \cite{M2}. Existence of such
a right parametrix with compact remainder implies that the range of $\Lbe$
is closed in all cases, and that its cokernel is at most finite dimensional
for a certain range of values of $\nu$. To show that this cokernel is trivial,
and that the norm of $\Gbe$ is uniformly bounded for sufficiently small
$\e_i$, we employ scaling arguments. These arguments proceed by contradiction,
showing that if first surjectivity, and secondly uniform surjectivity were to
fail, then counterexamples for increasingly small $\e_i$ could be
rescaled to obtain some element of the nullspace of the global operator
$ L_1$ on $\RNO$, or $\LL_1$ on $\Rnk$ which we show cannot exist.
Once these results are established, the rest of the proof of the
existence of solutions follows from a standard fixed point argument.

We remark here that we have proved these results concerning the uniform
surjectivity of $\Lbe$ in slightly greater generality than is used
in the nonlinear analysis. In particular, we are able to show that $\Lbe$
is uniformly surjective provided all $\e_i \le \e_0$ for $\e_0$ sufficiently
small. However, we are only able to prove the fixed point theorem if
$a\e < \e_i < \e$ for some $\e \le \e_0$, and for some $a \in (0,1)$,
i.e. when all the $\e_i$ are mutually relatively bounded.

\section{The globalized linearization}

In this section we analyze the behaviour of the operator
\[
L_1 = \Delta + pu_1^{p-1}
\]
on $\RNO$, and later the corresponding induced operator on
$\Rnk$.  In polar coordinates, $L_1$ takes the form
\begin{equation}
L_1 = \frac{\del^2}{\del r^2} + \frac{N-1}{r}\frac{\del}{\del r}
+ \frac{1}{r^2} \Delta_{\theta} + \frac{V_p(r)}{r^2},
\label{eq:L1-polar}
\end{equation}
where $V_p(r) = r^2 \cdot pu_1(r)^{p-1}$. From Proposition~\ref{pr:base} above
we have
\begin{eqnarray}
&\lim_{r \rightarrow 0} V_p(r) =  pv_{\infty}^{p-1} = \frac{2p}{p-1}
\left(N - \frac{2p}{p-1}\right) \equiv A_p, \\
&V_p(r) \sim cr^{(2-N)(p-1) + 2} \quad \mbox{as}\quad r \rightarrow \infty.
\end{eqnarray}
Notice that the exponent $(2-N)(p-1) + 2$ is negative precisely
when $p > N/(N-2)$.

\subsection{Indicial roots}
$L_1$ has a regular singularity at $r = 0$ and, because $V_p(r)$ tends
to $0$ as $r$ tends to infinity, it also has one at
$r = \infty$. Hence the asymptotic behaviour of solutions to
$L_1 w = 0$ are determined by the indicial roots of this
operator at these points. At $r=0$, the indicial roots for $L_1$ are
\begin{equation}
\gamma_j^{\pm} = \frac{2-N}2 \pm \sqrt {\left( \frac{N-2}{2} \right)^2
+ \lambda_j - A_p}.
\label{eq:indrts-0}
\end{equation}
Here the $\lambda_j$ are the eigenvalues of the Laplacian on
$S^{N-1}$, counted with their multiplicity. More precisely
 \[
\lambda_0 = 0, \quad \lambda_j= N-1,\ j = 1, \ldots N,\quad \mbox{etc.}
\]
The indicial roots for $L_1$
at $r = \infty$ are the same as for the Laplacian itself, since
$V$ tends to $0$ at $\infty$; these values are
\begin{equation}
\tilde \gamma_j^{\pm} = \frac{2-N}{2} \pm \sqrt
{\left( \frac{N-2}{2} \right)^2
+ \lambda_j}, \qquad j = 0, 1, 2, \dots.
\label{eq:indrts-inf}
\end{equation}
Notice that these numbers are integral and assume all integer values except
$-1, -2, \dots, 3-N$.

Because $L_1$ has a regular singularity at $0$ and $\infty$,
its mapping properties are well-known, cf.\ \cite{M}, \cite{Mc}.
\begin{proposition} The bounded linear map
\[
L_1: \cal C^{2,\alpha}_{\gamma,\gamma'}(\RNO)
\longrightarrow \cal C^{0,\alpha}_{\gamma-2,\gamma'-2}(\RNO)
\]
is Fredholm provided $\gamma \notin \{\gamma_j^{\pm}\}$ and
$\gamma' \notin \{\tilde \gamma_{j}^{\pm}\}$.
\end{proposition}

We will let $\gamma' = -1$ or $0$ since we are only interested in solutions
which decay at $\infty$, or at least are bounded. We could replace this
$-1$ by any value of $\gamma'$ in $(2-N,0)$. However, when
$N=3$, we would need to choose $\gamma' \in (-1,0)$. We shall not comment
on this further, and persist in letting $\gamma' = -1$, with this
understanding.

\subsection{Numerology}
We record some facts about these indicial roots and some of
the other constants that arise frequently for later reference.
\begin{lemma} Let $p$ be any number in the range $N/(N-2) <
p < (N+2)/(N-2)$.
\begin{enumerate}
\item The functions $-2/(p-1)$ and $A_p$ are monotone in $p$ and
\[
2-N > \frac{-2}{p-1} > \frac{2-N}2 \mbox{\qquad and \qquad}
0 < A_p \equiv \frac{2p}{p-1}\left(N - \frac{2p}{p-1}\right)
< \frac{N^2 - 4}{4}.
\]
\item There is a number $p^*$ in $(\frac{N}{N-2},\frac{N+2}{N-2}$ such that
for $N/(N-2) < p \le p^*$, the indicial roots $\gamma_0^{\pm}$ are
real, with $-2/(p-1) < \gamma_0^- < (2-N)/2 < \gamma_0^+$, while if
$p^* < p < (N+2)/(N-2)$, then $\gamma^{\pm}_0$ are both
complex with real part $(2-N)/2$.
\item $\gamma_1^- < -\frac{2}{p-1}$, in fact $\gamma_1^- =  -\frac{2}{p-1}
- 1$. Hence in particular, $\gamma_1^{\pm}$ are always real.
\end{enumerate}
\label{le:numerology}
\end{lemma}

\subsection{Injectivity of $L_1$ on $\C^{2,\alpha}_{\mu, 0}(\RNO)$}

Fix two weights $\nu$ and $\mu$ with
\begin{equation}
\frac{-2}{p-1} < \nu < \Re (\gamma_0^-)  \le \frac{2-N}2
\le \Re (\gamma_0^+)  < \mu
\label{eq:weights}
\end{equation}
and $\mu + \nu = 2-N$. These weights are `dual' in a sense to be
explained later. The interval in which $\mu$ lies is chosen to guarantee
injectivity of $L_1$ on $\C^{2,\alpha}_{\mu,0}$; $\nu$ is determined
by this duality, but also chosen not to be too negative in accordance
with the nonlinear problem.

\begin{proposition} The only solution $w \in \C^{2,\alpha}_{\mu,0}
(\RNO)$ of the equation $L_1 w = 0$ is the trivial solution $w = 0$.
\label{pr:is-inj-global}
\end{proposition}
\smallskip
{\bf Proof:} Let $\{\phi_j(\theta), \lambda_j\}$ be the eigenfunctions
and eigenvalues for $\Delta_{\theta}$, the Laplacian on $S^{N-1}$.
Then if $w$ is any solution of $L_1 w = 0$, it decomposes into
the infinite sum
\[
w(r,\theta) = \sum_{j = 1}^{\infty} w_j(r) \phi_j(\theta), \qquad
\]
where the $w_j$ solve
\[
L_{1,j} w_j \equiv
\frac{\del^2}{\del r^2} w_j + \frac{N-1}{r} \frac{\del}{\del r} w_j
+ \frac{1}{r^2}(V_p(r) - \lambda_j) w_j = 0.
\]
The value $j = 0$ is omitted in the sum above since every solution
of $L_{1,0}w_0 = 0$ blows up faster than $r^\mu$ near
$r = 0$.

Now, by the growth restrictions on $w$ near $0$ and infinity,
each $w_j$ must be bounded by $C_j r^{\gamma_j^+}$ as $r
\rightarrow 0$, and must decay like $r^{2-N-j}$ as $r \rightarrow
\infty$. We shall show that this forces $w_j$ to vanish.

First of all, notice that if $\lambda_j > \sup V_p$ then  solutions
of $L_{1,j}w_j = 0$ satisfy the maximum principle. For these
values of $j$ it is also true that $\gamma_j^+ > 0$, and so $w_j$ tends
to $0$ both at $0$ and infinity, hence must vanish identically.
Thus there are only finitely many values of $j$ for which $w_j$ may
possibly not vanish. To deal with most of these, we use an integral
estimate which will show that $w_j = 0$ for $\lambda_j > N-1$.
Multiply the equation $L_{1,j}w_j = 0$ by
$r^{N-1}w$ and integrate from $\kappa$ to $R$.  Integrate by parts
and use that $\mu > (2-N)/2$ and that $w_j$ decays at least as fast
as $r^{1-N}$ as $r \rightarrow \infty$ to see that the boundary
terms converge to zero as $\kappa \rightarrow 0$ and $R \rightarrow \infty$,
and that the integral converges. The result is that
\[
\int_0^{\infty} -(\del_r w_j)^2 r^{N-1} +
(V_p(r) - \lambda_j)w_j^2 r^{N-3}\,dr = 0.
\]
Now rearrange terms and estimate $V$ by its supremum:
\begin{eqnarray*}
\int_0^{\infty} r^{N-1} |\del_r w_j|^2\,dr & =
\displaystyle{\int_0^{\infty}} (V_p(r) - \lambda_j) r^{N-3}|w_j|^2\,dr \\
&\le (\sup V_p  - \lambda_j) \displaystyle{\int_0^{\infty}} r^{N-3} w_j^2\,dr.
\end{eqnarray*}
To proceed further, we note that there is a general inequality
for any function which decays like $r^{\mu}$ at $0$ and like $r^{1-N}$ at
infinity:
\[
\int_0^{\infty} r^{N-3}|w|^2\,dr \le \frac{4}{(N-2)^2} \int_0^{\infty}
r^{N-1}|\del_r w|^2\,dr .
\]
To prove this, observe that
\begin{eqnarray*}
 \int_0^{\infty} r^{N-3} w^2 \,dr =\frac{1}{N-2} \int_0^{\infty}  w^2
\del_r r^{N-2} \, dr = - \frac{2}{N-2} \int_0^{\infty} r^{N-2}w \del_r w
\, dr,
\end{eqnarray*}
where again the boundary terms vanish in the integration by parts
because of the assumed rates of decay.  Now use the
Cauchy Schwarz inequality to  prove the claim.

Use this inequality in the estimate for $w_j$ above to get
\[
\frac{(N-2)^2}{4}\int_0^{\infty} r^{N-3}  w_j^2\,dr
\le \left( \sup V_p - \lambda_j \right) \int_0^{\infty} r^{N-3}w_j^2\,dr.
\]
If $\sup V_p - \lambda_j < (N-2)^2/4$, then this inequality shows that
$w_j$ must vanish identically. By Proposition \ref{pr:base}) and Lemma
\ref{le:numerology}
\[
\sup V_p \le \frac{p+1}{2} \frac{2p}{p-1}\left( N - \frac{2p}{p-1}
\right).
\]
The function of $p$ on the right is dominated by its value at
$p = (N+2)/(N-2)$, so that
\[
\sup V_p < \frac{N}{N-2}\cdot \frac{N^2 - 4}{4} = \frac{N(N+2)}{4}.
\]
But now
\[
\frac{N(N+2)}{4} - \lambda_j < \frac{(N-2)^2}{4} \Longrightarrow
\lambda_j > \frac{3}{2}N - 1,
\]
and this occurs for all $j$ except $j = 0, 1, \ldots , N$.

Thus we have shown that $w_j = 0$ for $j > N$, and because of the
restrictions on the growth at $r = 0$, we also know that $w_0 = 0$.
It remains to show that $w_j = 0$, $j=1, \ldots N$.  It turns out that we can
explicitly write down an explicit solution of $L_{1,j}w = 0$ for
this range of $j$. In fact, differentiate the equation
$\Delta u_1 + u_1^p = 0$ with respect to $\del/ \del x_j$ to get
\[
L_1 \frac{\del u_1}{\del x_j} = 0.
\]
Since $u_1$ depends only on $r$, $\del u_1/\del x_j =
(\del u_1/\del r) \theta_j$, where $\theta_j = x_j / |x|$.
Since $-\Delta_{\theta} \theta_j = \lambda_j \theta_j$,
$j = 1, \ldots N$, $u_1'(r)$ solves $L_{1,j} w = 0$. But $u_1'(r)$ decays
like $r^{1-N}$ as $r \rightarrow \infty$ and blows up like $r^{-2/(p-1) - 1}$
as $r \rightarrow 0$. If $w_j$ solves $L_{1,j} w_j = 0$ and
decays like $r^{1-N}$ at infinity and like $r^{\gamma_j^+}$ at 0,
then some nontrivial linear combination of $u_1'$ and $w_j$
must decay faster than $r^{1-N}$ at infinity; since the singularities
of these functions at $0$ cannot cancel, this linear combination
is nonvanishing. This is a contradiction, since no solution of
$\cal L_{1,j}$ can decay faster than $r^{1-N}$ at infinity.
Hence $w_j = 0$, and the proof is complete. \hfill $\Box$

\subsection{Injectivity of $\LL_1$ on $\cal C^{2,\alpha}_{\mu,0}
(\Rnk)$}
The argument that the induced operator $\LL_1 = L_1 + \Delta_y$ on $\Rnk$
is injective on functions which blow up no faster than $r^{\mu}$
and decay like $r^{-1}$ as $r \rightarrow \infty$ is somewhat more
complicated, but rests on the results of the last section.

\begin{proposition} The only solution $w \in \cal C^{2,\alpha}_{\mu,0}
(\Rnk)$ which satisfies $\LL_1 w = 0$ is $w = 0$.
\label{prop:gen-inj-global}
\end{proposition}
\smallskip
{\bf Proof:} As in the previous subsection, we analyze this equation
by reducing it to a family of ODE's and studying these separately.
This is done by first taking the eigenfunction decomposition in $\theta$,
as before, but then also taking the Fourier transform in $y$.
Letting $\eta$ be the variable dual to $y$ in the Fourier transform,
we obtain the family of operators
\[
\LL_{j,|\eta|^2} = \frac{\del^2}{\del r^2} + \frac{N-1}{r}
\frac{\del}{\del r}
+ \frac{V_p(r) - \lambda_j}{r^2} - |\eta|^2.
\]
If $\LL_1 w = 0$ and $w \in \C^{2,\alpha}_{\mu,0}$, then
\[
w(r,\theta,y) = \sum_{j = 0}^{\infty} \int e^{iy\cdot \eta}\,
\hat{w}_j(r,\eta)\,dy,
\]
where $\LL_{j,|\eta|^2} \hat{w}_j = 0$ and $\hat{w}_j \in
\C^{2,\alpha}_{\mu,0}(\RR^+)$ and depends distributionally on $\eta$.
Clearly, if $\LL_{j,|\eta|^2}\hat{w}_j = 0$ with $\hat{w}_j
\in \C^{2,\alpha}_{\mu,0}(\RR^+)$ implies $\hat{w}_j = 0$ for every $j$,
then $\LL_1w = 0$ with $w \in \C^{2,\alpha}_{\mu,0}(\Rnk)$ implies
$w = 0$ as well.

First observe that the integration by parts argument of the last
subsection may be repeated essentially verbatim to show that
all solutions of this equation are trivial unless $j = 0, 1, \dots, N$.
Furthermore, since $\mu > \Re(\gamma_0^{\pm})$, there are no
solutions in $\C^{2,\alpha}_{\mu,0}(\RR^+ )$ for $j = 0$, either.
Thus we are left to consider only the case $j = 1, \ldots N$, where
$\lambda_j = N-1$.

For simplicity, in this argument, let us drop the index $j$ in our
notation and replace $|\eta|^2$ by a parameter $E$. Then $\LL_{j,|\eta|^2}$
becomes $\LL_{E}$, and so forth. Now let us recall that any solution of
$\LL_{E}w = 0$ has a Frobenius series expansion (convergent in some interval,
since $V_p$ is real analytic) of the form
\begin{equation}
w \sim \sum_{j = 0}^{\infty} a_j r^{\gamma^- + j} +
\sum_{j=0}^{\infty} b_j r^{\gamma^+ + j}.
\label{eq:frob-series}
\end{equation}
There are additional terms, with logarithmic factors, in the
second series in the case when the two indicial roots differ by
an integer. Since the argument is quite analogous in that case,
we only consider the case when these extra terms do not appear.
The two leading coefficients, $a_0$ and $b_0$ are free, and may
be specified arbitrarily; all the $a_j$ for $j > 0$ are determined
in terms of $a_0$ by a recursion formula, and similarly for the
$b_j$.  In particular, if $a_0 = 0$, then all $a_j = 0$. Hence,
since $\mu > \gamma_1^-$, we have that $w \in \cal C^{2,\alpha}_{\mu}$
near $r=0$ if and only if $a_0 = 0$.

We shall show by a sort of continuity argument, treating both $p$ and $E$ as
parameters, that if $w$ is a nontrivial solution on $\RR^+$
growing at most polynomially as $r \rightarrow \infty$ then $a_0$ never
vanishes.

For any two numbers $0 \le r_0 < r_1 < \infty$ there is a map
$C_{r_0,r_1}$ which sends the Cauchy data $(w(r_0),w'(r_0))$ of
a solution $w$ at $r_0$ to its Cauchy data $(w(r_1),w'(r_1))$
at $r_1$. When $r_0 = 0$, this Cauchy data is replaced by the
pair of leading coefficients $(a_0,b_0)$ in (\ref{eq:frob-series}).
If $0 < r_0$ then it is well-known that $C_{r_0,r_1}$ is
invertible and depends smoothly on the parameters $(p,E)$.
It is also true that if $r_1$ is in the interval of convergence of
the series (\ref{eq:frob-series}) then $C_{0,r_1}$ is invertible
and still depends smoothly on $(p,E)$.

Amongst the two dimensional family of solutions of $\LL_{E} w = 0$,
there is one solution, unique up to multiplication, which decays
at infinity. Actually, if $\LL_{E} w = 0$ and if $w$ is of polynomial
growth as $r \rightarrow \infty$, then $w$ decays like $r^q e^{-r\sqrt{E}}$
for $q$ determined by the coefficients of $\LL_{E}$.
We fix the solution uniquely by requiring that $w(1) = 1$.
This makes sense because, by the normalization in
\ref{re:small-outside}, the function $V_p(r)$ is less than some fixed
constant, in particular less than $N-1$, for $r \ge 1$. Thus solutions
of $\LL_{E} w = 0$ satisfy the maximum principle on $[1,\infty)$.
Hence any solution which vanishes at $r = 1$ cannot decay as
$r \rightarrow \infty$.

For this decaying solution $w$ look at the leading asymptotic coefficient
$a_0 = a_0(p,E)$ in its series expansion at $r=0$.
We write out this dependence on the two parameters $p$ and $E$ explicitly
since we need to show that $a_0(p,E) \ne 0$
for $(p,E)$ in the strip $(N/(N-2),(N+2)/(N-2)) \times [0,\infty)$.
For if this holds, then $w \notin \cal C^{2,\alpha}_{\mu,0}(\RR^+ )$ and our
result will follow.

We first show that $a_0$ is smooth in $(p,E)$. We do this by using
the maps $C_{r_0,r_1}$ discussed above. First consider the Cauchy
data $(1,w'(1))$ of the decaying solution $w$. Since the coefficient
of order zero of $\LL_{E}$ is negative on $[1,\infty)$, $w$ may be
constructed by the `shooting method' there, i.e. by solving the boundary
problem $\LL_{E} w_R = 0$ on $[1,R]$, with
$w_R(1) = 1$ and $w_R(R) = 0$, and letting $R \rightarrow \infty$.
By the maximum principle, the limit must exist, and since it is bounded
it agrees with $w$. Since the Cauchy data $(1,w_R'(1))$ depends smoothly
on $(p,E)$ for each $R$, it is not difficult to see that the limiting
values $(1,w'(1))$ do as well. Now use the composition of the (invertible)
maps $C_{0,r_0}$ for $r_0$ very small, and $C_{r_0,1}$ to conclude that
$(a_0,b_0)$ depends smoothly on $(p,E)$.

Now we proceed with the rest of the argument. First observe that
$a_0(p,0) \ne 0$ for any $p$ by the result of the previous subsection.
In fact, writing out the explicit solution
$\del u_1/\del r$ more explicitly, we can check that
$a_0(p,0) > 0$. Next, note that for $p$ sufficiently close
to its lower limit $N/(N-2)$, the function $V$ stays uniformly
small, and in particular less than $N-1$ for all $r$. Hence we
can again apply the integration by parts argument of the previous
subsection, but now using that $\sup(V_p - N+1) < 0$, which
immediately shows that $a_0 \ne 0$ for all $p < N/(N-2) + \kappa$
for every $E$ and for some $\kappa > 0$.
By continuity from $E = 0$, $a_0 > 0$ for
$(p,E)$ in this smaller strip. Finally, we can apply a similar argument
when $E$ is sufficiently large, regardless of the value of $p$.
In fact, all we need is that the supremum of $V_p(r) - N+1 - r^2 E$
is less than $(N-2)^2/4$, and this is clearly true for any $p$, for
$E$ large enough.  Thus, $a_0 > 0$ also for $(p,E)$ in this range
of values.

Finally, suppose that there is some value of $(p,E)$, with
$p \in (N/(N-2), (N+2)/(N-2))$, for which $a_0 = 0$. Using the
regularity of $a_0$ in these parameters, there is some $(p_0,E_0)$ for which
$a_0(p_0,E_0) = 0$ and $(a_0)_E (p_0,E_0) = 0$ (i.e. both $a_0$ and
its derivative with respect to $E$ vanish there).
Differentiate the equation $\LL_{E}w = 0$ and set $(p,E) = (p_0,E_0)$.
If $\tilde{w} = \del w / \del E$ then
\[
\LL_{E} \tilde{w} = 2w.
\]
Since $a_0 = 0$, the right side of this equation blows up only
like $r^{\gamma^+}$ as $r \rightarrow 0$. Thus $\tilde{w}$ is
the sum of two terms; the first blows up no faster than
$r^{\gamma^+ + 2}$, while the second is a solution of the homogeneous
equation and might blow up like $r^{\gamma^-}$. However, since the
leading coefficient of $\tilde{w}$ is the derivative of the leading
coefficient of $w$, i.e. it is $(a_0)_E$, which
vanishes, we see that this solution of the homogeneous
equation is absent and $\tilde{w}$ blows up no faster than
$r^{\gamma^+ + 2}$.  Now multiply the equation by $w \cdot r^{N-1}$
and integrate from $0$ to infinity. Using the exponential decay as $r$
tends to infinity, and that $w \sim r^{\gamma^+}$, $\tilde{w} \sim
r^{\gamma^+ + 2}$ as $r$ tends to $0$, we see that this integral is
well-defined, since $2\gamma^+ > 2-N$. Furthermore, it is permissible to
integrate by parts to get
\[
0 = \int_0^{\infty} (\LL_{E} w) \tilde{w}\, r^{N-1}\,dr =
\int_0^{\infty} w (\LL_{E} \tilde{w}) r^{N-1}\,dr =
2 \int_0^{\infty} |w|^2 r^{N-1}\,dr.
\]
Hence $w = 0$ for this value of $(p,E)$, which is a contradiction.
\hfill $\Box$

\subsection{Surjectivity of $L_1$ and $\LL_1$ on $\C^{2,\alpha}_{\nu,-1}$}
{}From the results of the last subsection we may deduce the following:
\begin{proposition}
Let $-\frac{2}{p-1} < \nu = 2-N - \mu < \Re (\gamma_0^-)$ as before. Then the
maps
\[
L_1: \C^{2,\alpha}_{\nu,-1}(\RNO) \longrightarrow \C^{0,\alpha}_{\nu-2,-3}
(\RNO)
\]
and
\[
\LL_1: \C^{2,\alpha}_{\nu,-1}(\Rnk) \longrightarrow \C^{0,\alpha}_{\nu-2,-3}
(\Rnk)
\]
are surjective. Furthermore, the first map is Fredholm and its nullspace
is one dimensional, while the second map has an infinite dimensional
nullspace.
\label{pr:surj-glob}
\end{proposition}
\smallskip
{\bf Proof:} We shall only sketch the proof of this result and provide
references to papers where analogous facts are proved thoroughly.
As discussed already in \S 5, the proof of this result
is much simpler for the former of the two maps.  In fact it is not difficult
to write down a right inverse for $L_1$ directly in terms of solutions
of the homogeneous problem $L_1 w = 0$. In order to do this, it is
crucial to use the fact that any solution which is bounded by a
multiple of $r^{\mu}$ near $r=0$ does not decay as $r$ tends to infinity,
and any solution which decays at infinity must decay at least as fast
as $r^{2-N}$. Note that this is equivalent to the injectivity of $L_1$ on
$\C^{2,\alpha}_{\mu,-1}$. The right inverse for $L_1$ may be written
as a sum of the right
inverses for $L_{1,j}$ on each eigencomponent; each of these right inverses
may be constructed explicitly from the solution which blows up like
$r^{\gamma^+_j}$ as $r \rightarrow 0$ which is unbounded as
$r \rightarrow \infty$ and the solution blowing up like $r^{\gamma^-_j}$
as $r \rightarrow  0$ and decaying as $r \rightarrow \infty$.
It is then not difficult to show that this sum of right inverses for each
eigencomponent is bounded on the weighted H\"older spaces,
cf. \cite{C-H-S} for this procedure for a closely related operator.
This proves that $L_1$ on $\C^{2,\alpha}_{\nu,-1}$ is not only
Fredholm but surjective.

Any solution of $L_1 w = 0$ in $\C^{2,\alpha}_{\nu,
-1}$ must be radial, by considering its growth as $r \rightarrow 0$.
There is a two dimensional space of solutions of $L_{1,0}w_0 = 0$.
By considering a suitable nontrivial linear combination of any two
basis elements of this space we obtain a solution which decays like
$r^{2-N}$ as $r \rightarrow \infty$; there cannot be two independent
solutions with this decay rate for the same reasons as in the previous
subsection. This solution grows like some combination of $r^{\gamma_0^{\pm}}$
as $r \rightarrow 0$, but in either case is in $\C^{2,\alpha}_{\nu,-1}$.
Hence the nullspace is one dimensional.

To establish the corresponding facts for $\LL_1$ we use Hilbert space
techniques and the construction of a pseudodifferential right parametrix
for $\LL_1$ from \cite{M2}. $\LL_1$ is surjective on
$\C^{2,\alpha}_{\nu,-1}$ provided its range is both closed and dense.
The existence and boundedness of this right parametrix  gives the
closedness of the range. Then duality, coupled with
Proposition~\ref{prop:gen-inj-global}, yields that the cokernel of this map is
trivial, so that the right parametrix may be replaced by a right inverse.

It is most natural to construct the pseudodifferential right parametrix
for $\LL_1$ relative to the spaces $r^{\delta}L^2$ because of the
central role of the Fourier transform in this construction.
Once the parametrix is obtained, it is then necessary to show its
boundedness on the weighted H\"older spaces. These steps are carried
out in detail in \cite{M2}, to which we refer the reader,
cf. also \cite{M-S}. Here we simply show
that $\LL_1$ satisfies the hypotheses necessary for that machinery to apply.

The main result of \cite{M2} implies that $\LL_1$ has closed range
on $r^{\delta}L^2$ provided $\delta \notin \{ \delta_j^{\pm}\}$,
where $\delta_j^{\pm} = \gamma_j^{\pm} + (N-2)/2$. Choose a
number $ -\delta_{\nu}$ just slightly smaller than
$\nu + (N-2)/2$, in particular so that
$-2/(p-1) + (N-2)/2 < -\delta_{\nu} < \nu + (N-2)/2 < 0$. Then
$-\delta_{\nu} \notin \{\delta_j^{\pm}\}$ and
$\C^{2,\alpha}_{\nu}(\OS) \subset \rho^{-\delta_{\nu}}L^2(\OS)$.
It is also proved in \cite{M2} that $\LL_1$ is
essentially surjective on $r^{-\delta_{\nu}}L^2$, i.e. it has
closed range and at most a finite dimensional cokernel, and that
it is essentially injective, i.e. has closed range and at most a
finite dimensional kernel, on $r^{+\delta_{\nu}}L^2$. These
two facts are consequences of one another by duality since $r^{\delta}L^2$
is the dual space of $r^{-\delta}L^2$ and $\LL_1$ is self-adjoint
on $r^{0}L^2$.  The crucial hypothesis that must be satisfied in
order for these last two conclusions
to be true involves the `model Bessel operator' $\hat{L}_1$ for $\LL_1$.
This operator is obtained by freezing coefficients of $\LL_1$ at $r=0$
in an appropriate sense, taking the Fourier transform in
$y$, and finally rescaling by setting $s = r|\eta|$. The operator
obtained in this way is
\[
\hat{L}_1 = \frac{\del^2}{\del s^2} + \frac{N-1}{s} \frac{\del }{\del s}
+ \frac{A_p + \Delta_{\theta}}{s^2} - 1.
\]
This crucial hypothesis is that $\hat{L}_1$ is injective as a map
on $s^{\delta_{\nu}}L^2(\RR^+ \times S^{N-1};s^{N-1}\,ds\,d\theta)$,
or equivalently, surjective as a map on $s^{-\delta_{\nu}}L^2(\RR^+
\times S^{N-1}; s^{N-1}\,ds\,d\theta)$.

$\hat{L}_1$ may be analyzed directly to show that this hypothesis is
satisfied. In fact, introducing the eigenfunction expansion with
respect to $\Delta_{\theta}$, as usual, we obtain a family of ordinary
differential operators, the solutions of which may be determined explicitly in
terms of Bessel functions. Analogously to the situation for $L_1$ we find
that any solution which blows up no faster than $r^{\nu}$ as $r \rightarrow
0$ grows exponentially as $r \rightarrow \infty$. Using these solutions
we can construct a right inverse for $\hat{L}_1$ explicitly.

Since $\hat{L}_1$ satisfies the hypothesis, we conclude that
$\LL_1$ itself has closed range, with at most a finite dimensional
cokernel, as a map on $\C^{2,\alpha}_{\nu,-1}$. The cokernel of
this map may be identified with the kernel of $\LL_1$ as a map
on $r^{\delta_{\nu}}L^2$ or $\C^{2,\alpha}_{\mu,-1}$. But from the
previous subsection, we know that this cokernel is trivial; hence
$\LL_1$ is surjective, as desired.

Finally, we may explicitly exhibit an infinite dimensional nullspace
of $\LL_1$ in $\C^{2,\alpha}_{\nu,-1}$.  For, if $\LL_1 w = 0$ and
$w$ is bounded by $r^{\nu}$ as $r \rightarrow 0$, then only
the eigencomponent $w_0$ may be nonvanishing, just as for $L_1$.
Now, taking the Fourier transform in $y$, we see that
$\LL_{0,|\eta|^2} \hat{w}_0 = 0$. Hence $\hat{w}_0(r,\eta) =
A(\eta) \hat{W}_0(r,\eta)$, where $\hat{W}_0$ is the unique solution
of this equation which
decays (exponentially if $\eta \ne 0$) as $r \rightarrow 0$,
normalized so that $\hat {W}_0(1,\eta) = 1$. The coefficient $A(\eta)$
is allowed to be arbitrary, so long as the corresponding $w_0$ is
bounded uniformly, for each $r$, as $|y| \rightarrow \infty$.
Clearly there is an infinite dimensional freedom in choosing such
coefficients; for example, we could let $A(\eta)$ be an arbitrary
element of the Schwartz space. This completes the proof.
\hfill $\Box$


\section{Injectivity of $\Lbe$ on $\cal C^{2,\alpha}_{\mu,\cal D}(\OS)$}

We are now able to turn to the first of our main tasks, to show that
$\Lbe$ is surjective on $\cal C^{2,\alpha}_{\nu,\cal D}(\Omega \setminus
\Sigma)$ when all the $\e_i$ are sufficiently small.
By an argument identical to the one indicated in the proof of
Proposition~\ref{pr:surj-glob}, this surjectivity is equivalent
to the injectivity of this operator acting on functions growing
like $r^{\mu}$ at the singular points, where $(2-N)/2 < \mu = 2- N -\nu <
\gamma_1^{+} <0$, as in the previous section. This injectivity is what we
shall prove.

\subsection{Preparatory Lemmas}
There are three lemmas which are used in the rescaling proofs in
both cases. We shall adopt the notation
\[
\Omega_{\be} = \Omega \setminus \cup_{i=1}^K B_{\e_i}(\Sigma_i).
\]
Here, using slightly different notation than in previous sections,
$B_{\e_i}(\Sigma_i)$ is the tubular neighbourhood around
$\Sigma_i$ of radius $\e_i$, and, in case $\Sigma$ is discrete,
$\Sigma_i$ simply equals the point $x_i$.

The first Lemma states that $\Lbe$ satisfies the maximum
principle on $\Obe$.
\begin{lemma} There exists some $\e_0 > 0$ such that if all $\e_i < \e_0$,
then after a normalization of the initial radial solution $u_1$,
the operator $\Lbe$ satisfies the maximum principle on
$\Obe$. That is, if $w \le 0$ on $\del \Omega_{\be}$ and $\Lbe w \ge 0$,
then $w \le 0$ on $\Omega_{\be}$.
\label{le:max-princ}
\end{lemma}
\smallskip
{\bf Proof:} We give the proof of this result in the case where $\Sigma$ is
discrete. The general case, where $\Sigma$ has positive dimension can be
treated similarly. Set $w^+ = \max\{w,0\}$ on $\Omega_{\be}$ and
$w^+ = 0$ on each $B_{\e_i}(\Sigma_i)$, so that $w^+ \in H^1_0(\Omega)$.
Multiply the inequality $\Lbe w \ge 0$ by $w^+$ and integrate
by parts to get
\[
\int_{\Omega} |\nabla w^+|^2 \le \int_{\Omega} |w^+|^2 p\bar{u}_{\be}^{p-1}
\le p\alpha^{p-1}\int_{\Omega}(\sum_{i=1}^K|x-x_{i}|^{-2}) |w^+|^2 .
\]
Here $\alpha$ is the supremum of $v_1 = r^{\frac{2}{p-1}}u_1$
for $r \ge 1$, which by Remark~\ref{re:small-outside} can be taken as
small as desired. Now, using the identity
\[
(N-2)|x|^{-2} = \mbox{div\ } (\frac{x}{|x|^2}),
\]
and integrating by parts, we obtain the inequality
\[
\int_{\Omega}|x-x_{i}|^{-2}|w^+|^2 \leq
\frac{4}{(N-2)^2} \int_{\Omega} |\nabla w^+|^2
\]
If $\alpha$ is taken sufficiently small so that the inequality
$4p\alpha^{p-1} K < (N-2)^2$ holds, we see that $w^+ = 0$, and hence $w
\le 0$, as desired.  \hfill $\Box$

\smallskip
The second Lemma uses this maximum principle to deduce estimates
in the weighted H\"older spaces for solutions of $\Lbe w_{\be} =
\fbe$ in terms of their boundary values on $\del \Obe$.

\begin{lemma}
There exists $\e_0>0$ such that, if all $\e_i < \e_0$ then the following
estimate holds. If $\fbe \in \C^{0,\alpha}_{\gamma - 2}(\Obe)$, where
$\gamma \in (2-N,0)$ is fixed, and suppose that $w_{\be}$ is  any solution
to $\Lbe w_{\be} = \fbe$ with $w_{\be} = 0$ on the `outer boundary'
of $\Obe$, i.e. on $\del \Obe \cap \del \Omega$.
Then there exists a constant $c>0$ independent of $\be$ such that
\[
||w_{\be}||_{2,\alpha,\gamma} \leq c \left(||\fbe||_{0,\alpha,\gamma-2}
+ \sum_{i=1}^K \e_i^{-\gamma}||w_{\be}||_{0,0,\del B_{\e_i}(\Sigma_i)}
\right),
\]
where the first two norms are taken over $\Omega_{\be}$.
\label{pr:est-outside}
\end{lemma}
\smallskip
{\bf Proof~:}
Define $\phi \in \C^{\infty}(\OS)$ to be a positive smooth function
for which, in some fixed neighbourhood $B_{\sigma}(\Sigma_i)$ for
each $i$, $\phi(x) = \mbox{dist\,}(x,\Sigma_i)^{\gamma}$; here $\sigma$
should also be chosen so that the approximate solution $\bue$ is supported
in $\cup B_{\sigma}(\Sigma_i)$. If some $\Sigma_i$ is a point, $\Sigma_i =
x_i$, then $\phi = |x-x_i|^{\gamma}$ in $B_{\sigma}(x_i)$. For example, we
can take $\phi = \rho^{\gamma}$, where $\rho$ was the function introduced
in \S3. Then a simple computation shows that in $B_{\sigma}(x_i) \setminus
B_{\e_i}(x_i)$ we have
\[
L_{\e_i}|x-x_i|^{\gamma} = \left\{\gamma (N-2+\gamma) |x-x_{i}|^{\gamma
-2}+p \bue^{p-1} |x-x_{i}|^{\gamma} \right\}
\leq - c |x-x_{i}|^{\gamma -2},
\]
where the constant $c > 0$ can be chosen independent of $\be$, since
$2-N < \gamma < 0$ implies that $\gamma(N-2+\gamma) < 0$ and
since $p\alpha^{p-1}$
can be chosen as small as desired. If $\Sigma_i$ is positive dimensional, then
a similar estimate holds. For this we simply need to use
the expression (\ref{eq:lap-fermi}) and the estimates for
$e_1$ and $e_2$ there.

Let $A>0$ denote the supremum of $w_{\be}$ on $\cup_{i=1}^{K} \del
B_{\sigma}(x_i)$. Hence $\phi$, multiplied by a suitable constant times
$\left(A + ||\fbe||_{0,\alpha,\gamma-2} + \sum_{i=1}^K \e_i^{-\gamma}
||w_{\be}||_{0,0,\del B_{\e_i}(\Sigma_i)}\right)$ is a supersolution for
the problem in $B_{\sigma}(x_i) \setminus B_{\e_i}(x_i)$. Likewise,
$-\phi$ multiplied by a similar constant is a subsolution.

We claim that, if $\e_0 >0$ is small enough, then $A$ is bounded
by a constant times
\[
\left(||\fbe||_{0,\alpha,\gamma-2}
+ \sum_{i=1}^K \e_i^{-\gamma}||w_{\be}||_{0,0,\del B_{\e_i}(\Sigma_i)}
\right).
\]
In order to prove this claim we will argue by contradiction
and assume that we have a sequence of counterexamples to this assertion.
Thus suppose that there is some sequence of $K$-tuples $\bel =
(\el_1, \dots ,\el_K)$ such that for some subset of the indices $1,
\dots , K$, the corresponding $\el_j$ tend to zero. For convenience, we
take this subset to be $\{1, \dots, J\}$, so that $\el_j \rightarrow 0$
for $j \le J$ and $\el_j \ge c > 0$ for $j > J$. In addition, up to a
subsequence we may assume that $\el_j$ converges to some $\e_j > c$ for all
$j > J$. Suppose also that for each such $\bel$ there is a function $w_{\bel}
\in \C^{2,\alpha}_{\gamma}(\Obel)$ with $\Lbe w_{\bel} = \fbel$ and
such that $w_{\bel} = 0$ on the `outer boundary' of $\Obe$ and $\fbel
\in \C^{0,\alpha}_{\gamma - 2}(\Obe)$ bounded independently of $\ell$.
First multiply the $w_{\bel}$ by a suitable constant so that
\begin{equation}
\sup_{\cup_{i=1}^{K}B_{\sigma}(x_i)} |w_{\bel}| = 1.
\label{eq:sup-omega=1}
\end{equation}
Using the super-solution constructed above, it is easy to see that
$w_{\bel}$ converges to a solution $w$ of $\Delta w =0$ in $\Omega
\setminus \cup_{i=1}^{J}\Sigma_i \cup_{i=J+1}^{K} B_{\e_i}(\Sigma_{i})$.
In addition, $w$ also vanishes on $\del \Omega\,\cup_{i=J+1}^{K}
\del B_{\e_i}(\Sigma_{i})$, and finally $w\in  \C^{2,\alpha}_{\gamma }$.
Since $\gamma > 2-N$ it is well known that the singularities of $w$ at
$\cup_{i=1}^{J} \Sigma_i $ are removable. Therefore $w$ solves
$\Delta w=0$ in $\Omega \setminus \cup_{i=J+1}^{K} B_{\e_i}(\Sigma_{i})$
with $0$ boundary data, so $w=0$ which contradicts
(\ref{eq:sup-omega=1}). The claim follows.

We conclude that there is a constant $\delta > 0$ independent of $\be$
such that $\Lbe (w_{\be} + \delta \phi) \leq f_{\be} - c\delta \phi \leq 0$
and $\Lbe (w_{\be} - \delta\phi)\geq \fbe +  c\delta\phi \geq 0$. By
increasing $\delta$ to make $w_{\be} -\delta\phi < 0$
and $w_{\be} + \delta\phi > 0$ on $\del \Obe$, we deduce that
$|w_{\be}| \leq \delta\phi$ in all of $\Obe$, i.e. $w_\e \in
{\C}^{0,0}_{\nu}(\Omega_{\be})$, with norm independent of $\be$.
Standard arguments using rescaled elliptic estimates show that $w_{\be}
\in {\C}^{2,\alpha}_{\nu} (\Omega_{\be})$ with norm independent of
$\be$. \hfill $\Box$

\smallskip

The third Lemma shows that the weighted H\"older norm with
exponent $\mu$ for a solution of the fixed homogeneous equation $L_1 w = 0$
on a fixed tubular neighbourhood $B_{\sigma}(\Sigma_i)$ is controlled
by the norm of $w$ on $\del B_{\sigma}(\Sigma_i)$. This is quite obvious
when $\Sigma_i$ is a point, since the nullspace of $L_1$ is finite
dimensional then.  When $\Sigma_i$ has positive dimension, this is
a more substantial result, and it follows directly from the results of
\cite{M2}. Nonetheless, we give an elementary proof, which covers
both cases. For
convenience here we drop the subscript $i$. As usual, $(r,\theta,y)$ are
Fermi coordinates on this tubular neighbourhood.

\begin{lemma} Suppose $w_{\ell}$ is a sequence of solutions of $L_1
w_{\ell} = 0$ in $B_{\sigma}(\Sigma)$, with $w_{\ell} \in \C^{2,
\alpha}_{\mu}(B_{\sigma}(\Sigma))$ and $|w_{\ell}| \le A$ on $\del
B_{\sigma}(\Sigma)$ uniformly in $\ell$.
Then $||w_{\ell}||_{2,\alpha,\mu}$ cannot diverge as $\ell \rightarrow
\infty$.
\label{le:bound-control}
\end{lemma}
\smallskip
{\bf Proof:} First of all, observe that by the rescaled Schauder
estimates, it suffices to show that the supremum of
$r^{-\mu}|w_{\ell}|$ over $B_{\sigma}$ does not diverge.
Suppose, to the contrary, that it does; furthermore, suppose
that this supremum takes the value $C_{\ell}$, which tends
to infinity, and is attained at some point $x_{\ell} \in B_{\sigma}$.
Let $(r_{\ell}, \theta_{\ell}, y_{\ell})$ be the Fermi coordinates
of $x_{\ell}$. Since $w_{\ell}$ can be bounded by some new  constant
$A'$ on $B_{\sigma}(\Sigma) \setminus B_{\sigma/2}(\Sigma)$, the
$r_{\ell}$ must converge to zero. For if they did not, then
consider the rescaled function $\bar{w}_{\ell} = C_{\ell}^{-1}w_{\ell}$.
This still solves $L_1 \bar{w}_{\ell}=0$. Furthermore, the supremum of
$r^{-\mu}\bar{w}_{\ell}$ on $B_{\sigma}(\Sigma)$
is equal to one, and is attained in some fixed annulus $B_{\sigma}(\Sigma)
\setminus B_{\beta}(\Sigma)$. Thus, the $\bar{w}_{\ell}$ converge to
a nonzero limit $\bar{w}_{\infty}$, which is a solution of
$L_1 \bar{w}_{\infty} = 0$ in $B_{\sigma}(\Sigma)$. However,
the supremum of $\bar{w}_{\ell}$ on $B_{\sigma}\setminus B_{\sigma/2}$
tends to zero, so that $\bar{w}_{\infty} = 0$ on this annulus,
which is a contradiction.

Now, since we have established that $r_{\ell} \rightarrow 0$, we may
again rescale. If $\Sigma = x_i$ is a point (which is the origin
in these coordinates), consider the function $\tilde{w}_{\ell}(r, \theta) =
C_{\ell}^{-1}r_{\ell}^{\mu}w_{\ell}(r/r_{\ell},\theta)$. If $\Sigma$
is higher dimensional, then some subsequence of the points $y_{\ell}$
converge to $y_{\infty} \in \Sigma$ (since $\Sigma$ is compact).
Choose Fermi coordinates centred around this point, so that
$y_{\infty} = 0$, and these coordinates are defined for $|y| \le
\tau$. In this case let $\tilde{w}_{\ell} = C_{\ell}^{-1}r_{\ell}^{\mu}
w_{\ell}(r/r_{\ell}, \theta, (y - y_{\ell})/r_{\ell})$. In the former
case, $\tilde{w}_{\ell}$ satisfies $L_{1/r_{\ell}} \tilde{w}_{\ell} = 0$
on $B_{\sigma/r_{\ell}}(0)$, and is bounded by $r^{\mu}$ there.
The same is true in the latter case, except that the operator $L_{1/r_{\ell}}$
must be replaced by one for which the error terms $e_1$, $e_2$
are replaced by very small translates, and which are in any case
still tending to zero when $\ell \rightarrow \infty$.
This rescaling has been chosen so that $r^{-\mu}|\tilde{w}_{\ell}|$
attains its supremum on $\del B_{1}(\Sigma)$. As before,
pass to a limit, $\tilde{w}_{\infty}$. By the previous remark,
$\tilde{w}_{\infty} \ne 0$.

When $\Sigma = 0$, then it is also true that $\tilde{w}_{\infty} \in
\C^{2,\alpha}_{\mu,\mu}(\RNO)$, and is a solution of
\[
\left(\Delta + \frac{1}{r^2} pv_{\infty}^{p-1}\right) \tilde{w}_{\infty} = 0
\]
there. But the solutions of this equation are of the form
\[
\sum_{j = 0}^{\infty} \left(a_j  r^{\gamma_j^+}
+ b_j r^{\gamma_j^-}\right)\phi_j(\theta).
\]
It is clear that no function of this form can be bounded by
$r^{\mu}$ both as $r \rightarrow 0$ and $r \rightarrow \infty$.
Hence we arrive at a contradiction again, so the assertion
of the lemma must be true when $\Sigma$ is a point.

In case $\mbox{dim\,}\Sigma = k > 0$, we may still take a limit,
and get a function $\tilde{w}_{\infty}$ which solves the same equation,
but now on ${\RR}^n \setminus {\RR}^k$. It also satisfies the
estimate $\sup r^{-\mu}|\tilde{w}_{\infty}| < \infty$. However,
introducing a decomposition into eigenfunctions for the spherical
Laplacian $\Lap_{\theta}$ yields the uncoupled system of equations
\[
\left(\frac{\del^2}{\del r^2} + \frac{N-1}{r} + \Lap_{{\RR}^k} +
\left(pv_{\infty}^{p-1} - \lambda_j\right)\frac{1}{r^2}\right)w_j = 0,
\]
where we have renamed the eigencomponents simply $w_j$. Each
such eigencomponent satisfies $\sup r^{-\mu}|w_j| < \infty$
uniformly as $|y| \rightarrow \infty$. Taking the Fourier transform of this
equation in $y$ reduces it to an equation of Bessel type:
\[
\left(\frac{\del^2}{\del r^2} + \frac{N-1}{r}\frac{\del}{\del r} +
\left( A_p - \lambda_j - r^2|\eta|^2\right) \frac{1}{r^2} \right)\hat w_j = 0.
\]
This equation may be solved explicitly, and as discussed earlier
it is easy to see that for $\eta \ne 0$
the only solution of this equation
which grows no faster than $r^{\mu}$ as $r \rightarrow 0$
must increase exponentially as $r \rightarrow \infty$.
Thus $\hat w_j(r,\eta) = 0$ for $\eta \ne 0$, and so it is polynomial,
and hence constant in $y$ (since $w_j$ is bounded in $y$). Therefore,
the previous argument applies to show $w_j = 0$ for all $j$.
Hence $\tilde{w}_{\infty} = 0$, which is again a contradiction.
This proves the Lemma in all cases. \hfill $\Box$
\smallskip

\subsection{Injectivity}

\begin{theorem}  If $\e_0$ is sufficiently
small, and if each  $\e_i < \e_0$, then
\[
\Lbe: \C^{2,\alpha}_{\mu,\cal D}(\OS)
\longrightarrow \C^{0,\alpha}_{\mu - 2}(\OS)
\]
is injective.
\label{th:inj-omega}
\end{theorem}
\smallskip
{\bf Proof:}  As in the previous lemma, the proof uses rescaling to argue
by contradiction. The idea is quite simple. If the result of this theorem
were false, there would exist a sequence of counterexamples, corresponding
to some sequence of the $\bel$, with at least some subset of the $\el_j$
decreasing to zero.  By rescaling these functions and passing to a limit,
a counterexample to Propositions~\ref{pr:is-inj-global} and
\ref{prop:gen-inj-global} would be obtained.

Thus suppose that there is some sequence of $K$-tuples $\bel =
(\el_1, \dots \el_K)$ such that (possibly passing to a subsequence)
for some fixed subset of the indices $1,\dots , K$, the corresponding $\el_j$
tend to zero. For convenience as before, we take this subset to be $\{1,
\dots, J\}$, so that $\el_j \rightarrow 0$ for $j \le J$ and $\el_j \ge c >
0$ for $j > J$. Suppose also that for each such $\bel$ there is a function
$w_{\ell} \in \C^{2,\alpha}_{\mu,\cal D} (\OS)$ with $\Lbel w_{\ell} = 0$.
First multiply the $w_{\ell}$ by a suitable constant so that
\begin{equation}
\sup_{\del \Obel} \rho(x)^{-\mu}|w_{\ell}| = 1.
\label{eq:sup-omega}
\end{equation}
 From Lemma~\ref{pr:est-outside}, we get that
\[
\sup_{\Obel} \rho(x)^{-\mu}|w_{\ell}| \le C\sup_{\cup \del B_{\el_i}(
\Sigma_i)}\rho^{-\mu}|w_{\ell}| \le C.
\]
For some subsequence (which we assume is the full sequence), the supremum
of $\rho(x)^{-\mu}|w_{\ell}|$ on $\cup \del B_{\el_i}(\Sigma_i)$ is attained
on some fixed $\del B_{\el_j}(\Sigma_j)$.

There are several possibilities, which we examine in turn. For simplicity,
we initially consider the case where $\Sig$ is discrete.
First suppose that for the sphere where supremum above is attained,
the index $j$ satisfies $j \le J$. Fix this index $j$, and for convenience,
translate $\Omega$ so that $x_j = 0$. Then $\el_j \rightarrow 0$.
Now rescale, setting
\[
\bar{w}_{\ell}(x) = (\el_j)^{-\mu}\,w_{\ell}(\el_jx).
\]
This function is defined on $B_{\sigma/\el_j}(0)$ and solves $(\Delta +
p\bar{u}_1^{p-1})\bar{w}_{\ell} = 0$ there. Here $\bar{u}_1$ is the radial
solution $u_1(r)$ multiplied by a cutoff function $\chi(\el_j x/R)$, which
is equal to one on an increasingly large set as $\el_j$ tends to
zero. By construction, the supremum of $r^{-\mu}|\bar{w}_{\ell}|$
is bounded by $C$ on $B_{\sigma/\el_j}(0) \setminus B_1(0)$.

We can apply Lemma~\ref{le:bound-control} to conclude that the
$r^{-\mu}|\bar{w}_{\ell}|$ are bounded uniformly on $B_1(0)$. Now let
$\ell \rightarrow \infty$. Since their norms are bounded, the
$\bar{w}_{\ell}$ tend to a limit $\bar{w}_{\infty} \in \C^{2,\alpha}_{\mu}
(B_1(0))$, which is defined with $r^{-\mu}|\bar{w}_{\infty}|$ bounded
on all of $\RNO$. In addition, it is a solution of $L_1\bar{w}_{\infty}= 0$
there. But $\bar{w}_{\infty}$ is not identically zero,
since its supremum on $\del B_{1}(0)$ is one. This contradicts
Proposition~\ref{pr:is-inj-global}.

This last argument still leads to a contradiction if it is only true that
for some $j \le J$ (and some subsequence of the $\ell$), the supremum of
$\rho^{-\mu}|w_{\ell}|$ on $\del B_{\el_j}(x_j)$ is bounded below by
some positive number $C' > 0$, independently of $\ell$.

Having excluded these cases, we can now assume that
\[
\sup_{\del B_{\el_j}(x_j)} \rho(x)^{-\mu}|w_{\ell}| \longrightarrow 0
\]
every $j \le J$.

The norms $||w_{\ell}||_{2,\alpha,\mu}$ (over all of $\Omega$) are bounded.
This is clear in the region $\Omega_{\bel}$ by (\ref{eq:sup-omega}), and
by Lemma~\ref{le:bound-control} it is also true for the balls
with index $j > J$. For the balls with index $j \le J$, note that the
rescaled functions $\bar{w}_{\ell} = (\el_j)^{-\mu}w_{\ell}(\el_j x + x_j)$
on $B_{1}(0)$ have the same weighted norm there as the unrescaled
functions do over $B_{\el_j}(x_j)$. Since the supremum of
$\bar{w}_{\ell}$  over $\del B_{1}(0)$ is bounded by $(\el_j)^{-\mu}$, which
tends to zero as $\ell \rightarrow \infty$ (and in particular, is
bounded), we can apply Lemma \ref{le:bound-control} again to conclude
boundedness of the weighted norms of $w_{\ell}$ in these balls too.

Now (pass to an appropriate subsequence and) let $\ell$ tend to infinity.
The $w_{\ell}$ converge to a nonzero limit $w_{\infty}$, since they have
supremum equal to one on $\del B_{\el_j}(x_j)$ for some $j > J$,
i.e. on a sphere which has radius bounded away from zero. This limit is an
element of $\C^{2,\alpha}_{\mu,\cal D}(\OS)$, and solves
$L_{\be'} w_{\infty} = 0$, where $\be'$ has $\e_j' = 0$ for $j \le J$.
Hence the potential $p\bar{u}^{p-1}_{\bar{\e}'}$ in $L_{\be'}$ is singular
only at the points $x_j$ for $j > J$.
Since $w_{\infty} \in \C^{2,\alpha}_{\mu}$ locally near each $x_j$,
and since $\mu > 2-N$, it easy to check that $w_{\infty}$ is a weak solution
of $L_{\be'}w = 0$ in a full neighbourhood of the points $x_j$ with
$j \le J$. Since the operator is smooth at these points, a standard removable
singularities theorem shows that $w_{\infty}$ is smooth except at the $x_j$,
$j > J$. This means that $w_{\infty}$ is a solution of this operator
with singularities at some discrete set $\Sigma'$ with
strictly fewer elements than $\Sigma$. Now we may proceed
by induction, the case $K=1$ already having been treated by
the proof above. This completes the proof when $\Sig$ is
discrete.

The modifications necessary to handle the general case are rather minor.
In fact, it is only necessary to modify the way in which the
rescaling is done and invoking Proposition~\ref{prop:gen-inj-global}
at the appropriate place. Thus, starting the proof in the same
way, if the functions $r^{-\mu}|w_{\ell}|$ attain their maximum
at a point $z_{\ell} = (r_{\ell},\theta_{\ell},y_{\ell})$ in Fermi
coordinates around some $\Sig_i$, then we use the rescalings
\[
\bar{w}_{\ell}(r,\theta,y) = (\e_j^{(\ell)})^{\mu} w_{\ell}(r/\el_j,\theta,
(y-y_{\ell})/\el_j).
\]
The rest of the proof is analogous to the previous case. \hfill $\Box$
\smallskip

\begin{remark}
It is actually possible to prove that in the case of nonisolated
singularities, the rescaled limiting function $\bar{w}_{\infty}$
is independent of $y$. This may be done by using some regularity
results from \cite{M} to show that the tangential oscillation of the
initial sequence $w_{\ell}$ may be bounded independently of $\be$.
\end{remark}

\section{Uniform surjectivity of $\Lbe$ on
$\C^{2,\alpha}_{\nu,\cal D}(\OS)$}

The second main step in the linear analysis is to show that there is
some choice of right inverse $\Gbe$ of $\Lbe$ on $\C^{0,\alpha}_{\nu-2}
(\OS)$ which has norm bounded independently of $\be$ provided each
$\e_j < \e_0$.  As we have indicated before, the subtlety here is
that $\Gbe$ is not unique, since $\Lbe$ is not injective on
$\C^{2,\alpha}_{\nu}$.

The usual choice for $\Gbe$ is as the right inverse whose range equals the
orthogonal complement of the nullspace of $\Lbe$. Of course, since we are
working in H\"older spaces, this orthogonal complement is meaningless.
However, what amounts to the same thing is to require that the range of
$\Gbe$ lies in the range of a fixed adjoint $\Lbe^*$ of $\Lbe$. Again, this
adjoint depends on some Hilbert space structure, but once we have chosen an
appropriate one with respect to which the adjoint is taken, we can forget
about it and simply use this adjoint.

By (\ref{eq:space-incl}) and using the notation of \S 6.5,
$\C^{2,\alpha}_{\nu}$ is contained in
$\rho^{-\delta_{\nu}}L^2(\rho^{N-1}dr\,d\theta\,dy)$.
In the following, the space $L^2$ will always be taken relative to the
background Euclidean measure. We shall consider the spaces
$\rho^{\delta}L^2$ and $\rho^{-\delta}L^2$ to be dual with respect
to the natural pairing
\[
\rho^{\delta}L^2 \times \rho^{-\delta}L^2 \ni (\phi,\psi) \longrightarrow
\int \phi \psi.
\]
Relative to this pairing, the adjoint of $\Lbe: \rho^{-\delta_{\nu}}L^2
\rightarrow \rho^{-\delta_{\nu} - 2}L^2$ is just $\Lbe: \rho^{\delta_{\nu}+2}
L^2 \rightarrow \rho^{\delta_{\nu}}L^2$. We have proved that the former
of these maps is surjective and the latter one is injective. Using the
fixed isomorphisms
\[
\rho^{2\delta}: \rho^{-\delta}L^2 \longrightarrow \rho^{\delta}L^2,
\]
we may identify this adjoint, $\Lbe^*$, with $\rho^{-2\delta_{\nu}} \Lbe
\rho^{2\delta_{\nu} + 4}$. Premultiplying by $\rho^{-4}$ we can regard
\[
\Lbe^* = \rho^{-2\delta_{\nu}}L \rho^{2\delta_{\nu}}:
\rho^{-\delta_{\nu} + 2}L^2 \longrightarrow \rho^{-\delta_{\nu}}L^2.
\]
It is not hard to see that Lemma~\ref{le:max-princ} is valid for both $\Lbe$
and $\Lbe^*$.

Now form the fourth order operator
\[
\cal L_{\be} = \Lbe \circ \Lbe^* = \Lbe \rho^{-2\delta_{\nu}}
\Lbe \rho^{2\delta_{\nu}}: \rho^{-\delta_{\nu}+2}L^2 \longrightarrow
\rho^{-\delta_{\nu} -2}L^2.
\]
This map is an isomorphism. Hence there exists a bounded two-sided inverse
\[
\cal G_{\be}: \rho^{-\delta_{\nu} - 2}L^2 \longrightarrow \rho^{-\delta_{\nu}
+ 2}L^2.
\]
In particular, $\cal L_{\be} \cal G_{\be} = I$, which means that
$\Gbe = \Lbe^* \cal G_{\be}$ is a right inverse to $\Lbe$ which maps
into the range of $\Lbe^*$ as desired.
Henceforth, $\Gbe$ will always denote this particular right inverse.

It can be shown directly from the structure of $\cal G_{\be}$ and
$\Gbe$ (see \cite{M2}) that
\[
\cal G_{\be}: \C^{0,\alpha}_{\nu-2}(\OS) \longrightarrow
\C^{4,\alpha}_{\nu + 2,\cal D}(\OS)
\]
and
\[
\Gbe: \C^{0,\alpha}_{\nu - 2}(\OS) \longrightarrow
\C^{2,\alpha}_{\nu,\cal D}(\OS)
\]
are bounded. In the former of these the subscript $\cal D$ denotes
a more general set of elliptic boundary conditions for this fourth order
operator, i.e. the domain is restricted to functions $u$ which vanish
along with $\Lbe^* u$ at $\del \Omega$.

We can now turn to the proof that the norm of $\Gbe$ does not blow up as the
components of $\be$ tend to zero. This uniform surjectivity is an immediate
consequence of the following two results.

\begin{proposition} Let $\be$ be some K-tuple with all $\e_j < \e_0$.
Then provided $\e_0$ is sufficiently small, there is no solution
of the system of equations $\Lbe u = 0$, $u = \Lbe^* v$ with
$u \in \C^{2,\alpha}_{\nu,\cal D}(\OS)$ and $v \in \C^{4,\alpha}_{\nu+2,
\cal D}(\OS)$.
\label{pr:inj-fourthorder}
\end{proposition}
\smallskip
{\bf Proof:}  This is really just a corollary of Theorem \ref{th:inj-omega}.
In fact, suppose that $u$ and $v$ satisfy this system. Then
$\Lbe \Lbe^* v = 0$. Recalling that $\Lbe^*$ is identified with
$\rho^{-2\delta_{\nu}} L \rho^{2\delta_{\nu}}$ here, multiply this
equation by $\rho^{2\delta_{\nu}}v = w$ and integrate with respect to
standard Euclidean measure. Then the integrations by parts in
\[
\int w \, \Lbe \rho^{-2\delta_{\nu}} \Lbe w = \int |\Lbe w|^2
\rho^{-2\delta_{\nu}} = 0
\]
are valid because $\del_{r}^j w \in \C^{4-j,\alpha}_{\nu + 2\delta_{\nu}
+ 2 - j}$ for $j \le 4$, and also because $w$ satisfies Dirichlet
conditions on $\del \Omega$. (Here $r$ is the Fermi polar distance coordinate
near each $\Sigma_i$. Also note that we can choose $\rho$ to be constant near
$\del \Omega$ so that $\rho^{2\delta_{\nu}}v$ and $\Lbe^*(\rho^{2
\delta_{\nu}}v)$ still vanish at the boundary.) Hence $\Lbe w = 0$.
But since $w \in \C^{2,
\alpha}_{\mu'}$ for some $\mu' > \Re (\gamma_0^{\pm})$, we may conclude from
Theorem \ref{th:inj-omega} that $w = 0$ provided $\e_0$ is small enough.
\hfill $\Box$
\smallskip

\begin{theorem} Let $\bel$ be any sequence of K-tuples, with
each $\el_i \le \e_0$. Let $\fbel$ be any sequence of functions
in $\C^{0,\alpha}_{\nu-2}(\OS)$ with norm uniformly bounded as $\ell$
tends to infinity. Let $\ubel \in \C^{2,\alpha}_{\nu,\cal D}
(\OS)$ be the unique solution of $\Lbe\ubel = \fbel$ which also satisfies
$\ubel \in \mbox{ran\ }(\Lbel^*)$, i.e. $\ubel= \Lbel^* \vbel$
for some $\vbel \in \C^{4,\alpha}_{\nu+2,\cal D}$. Then
the norm of $\ubel$ in $\cal C^{2,\alpha}_{\nu}$ is bounded
uniformly as $\ell \rightarrow \infty$.
\label{th:unif-surj}
\end{theorem}
\smallskip
{\bf Proof:} This argument is again by contradiction, and is
very similar to the ones in the last section. Clearly we only
need to consider the case where some subset of the $\el_j$ tend
to zero. As before, we assume that $\el_j \rightarrow 0$ for
$j = 1, \dots , J$ and $\el_j \ge c > 0$ for $j = J+1, \dots, K$.
To start, by hypothesis
\[
\sup_{\Omega} \rho^{-\nu +2} |\fbel| \le C,
\]
for all $\ell$. Applying Proposition~\ref{pr:est-outside}, we get
\[
\sup_{\Omega \setminus \cup B_{\el_i}(x_i)} \rho^{-\nu} |\ubel|
\le C + C \sup_{\cup \del B_{\el_i}(x_i)} \rho^{-\nu} |\ubel|.
\]
Now define
\[
A_{\be} = \sup_{\cup \del B_{\el_i}} \rho^{-\nu}|\ubel|.
\]
We resize the functions $\ubel,\ \vbel,\ \fbel$ by setting
\[
\tubel = \ubel/A_{\be},\quad \tvbel = \vbel / A_{\be}, \quad
\tfbel = \fbel / A_{\be},
\]
so that $\Lbel \tubel = \tfbel$ and $\Lbel^* \tvbel = \tubel$.
If $A_{\be}$ stays bounded as $\ell$ tends to infinity, then
we are finished. If not, as we now assume, then $||\tfbel||_{0,\alpha,\nu-2}$
tends to zero.

We wish to take a limit of the equations $\Lbel \tubel = \tfbel$
and $\Lbel^* \tvbel =\tubel$. The main point is to show that $\tubel$
and $\tvbel$ tend to limits. For each $\ell$, choose a point $z_{\ell}$ where
$\rho^{-\nu}\tubel$ attains its maximum. By passing to a subsequence,
we can assume that $z_{\ell}$ stays in some fixed $B_{\bel_i}(\Sigma_i)$.

The simpler case is when $i > J$, so that $\e_i$ does not tend to
zero. Since by the above and Lemma \ref{le:bound-control} the supremum
of $\rho^{-\nu}|\tubel|$ in $B_{\bel_i}(\Sigma_i)$ cannot diverge
and stays bounded away from zero.  Then we can directly pass to a limit
as $\ell$ tends to infinity and obtain a function $u \in \C^{2,\alpha}_{\nu,
\cal D}(\OS)$ which is nonvanishing and satisfies $L_{\bar{\kappa}}u = 0$.
Here $\bar{\kappa}$ is the $(K-J)$-tuple comprised of the limiting values
of the $\el_j$ for $j > J$. Since $u$ blows up no faster than $\rho^{\nu}$
at the $\Sigma_j$ with $j \le J$, and since the limiting term of order zero
in the operator is smooth at these submanifolds, the same removable
singularities theorem as we used before shows that $u$ must be smooth at those
points. Next, since the norm of $\tvbel$ is uniformly bounded, we can also
pass to a (weak) limit and obtain some function $v \in \C^{4,\alpha}_{\nu+2,
\cal D}(\OS)$ with $L_{\bar{\kappa}}^* v = u$. Since $u$ is nonvanishing, $v$
must also be. In addition, by elliptic regularity, $v$ is smooth at the
$\Sigma_j$, $j \le J$. However, this is a contradiction to Proposition
\ref{pr:inj-fourthorder}. Thus we may assume that for the index $i$ for which
the maximum of $\rho^{-\nu}|\tubel|$ is attained in $B_{\bel_i}(\Sigma_i)$, the
radius $\el_i$ tends to zero.

Let us first treat the case where $\Sigma_i$ is a point. Translate the
whole problem so that $\Sigma_i = \{0\}$. Next, rescale the two equations
by the factor $\bel_i$; in particular, $\tubel$ is replaced by
\[
\bubel = (\bel_i)^{-\nu}\tubel(\bel_i x),
\]
and similarly for $\tvbel$ and $\tfbel$. For convenience, replace $\bel_i$
by $\e$. Then in the ball $B_{\sigma/\e}(0)$ we have $\bfbel$ tending to
zero in $\cal C^{0,\alpha}_{\nu-2}$, and $\bubel$ having norm in
$\cal C^{2,\alpha}_{\nu}$ bounded uniformly by one. Replace
the rescaling of the function $\rho$ now by the polar variable $r$
in this expanding sequence of balls. In the new rescaled coordinates,
let the supremum of
$r^{-\nu}|\bubel|$ be attained at some point, which we
still call  $x_{\ell}$. Then $|x_{\ell}|$ is bounded. By the
argument used in the previous section, we see that $|x_{\ell}|$
must also stay bounded away from zero, otherwise we could rescale
by the factor $|x_{\ell}|^{-1}$ and arrive at a contradiction,
as in that section. Thus, there is some uniform lower bound
$0 < C \le |x_{\ell}| \le 1$, and hence we can pass to a limit
of the $\bubel$ in $\cal C^{2,\alpha}_{\nu}$.  Let us call the
limiting function $\bar{u}$. Then
\[
L_1 \bar{u} = 0 \quad \mbox{in}\quad \Bbb R^N \setminus \{0\}.
\]

We wish to show that $\bvbel$ also tends to a limit.  To see this,
first note, that by the maximum principle again,
\[
\sup_{B_{\sigma/\e}(0) \setminus B_1(0)} r^{-\nu} |\bvbel|
\le C + C \sup_{\del B_1(0)} |\bvbel|.
\]
Define
\[
A'_{\e} = \sup_{B_1(0)} r^{-\nu} |\bvbel|.
\]
If $A'_{\e}$ stays bounded as $\e \rightarrow 0$, we can
take a limit of the $\bvbel$ to get a function $v \in \C^{4,\alpha}_{\nu+2}
(\RNO)$ such that $L_1^* v = u$ there.

So, suppose not, i.e. suppose $A'_{\e}$ tends to infinity. We resize all the
functions once again, by letting $u_{\bel}' = \bubel/A'_{\e}$, and so on.
Then $L_1 u_{\bel}' = f_{\bel}'$ and $L_1^* v_{\bel}' = u_{\bel}'$.
By construction, the norm of $f_{\bel}'$ tends to zero in
$\C^{0,\alpha}_{\nu-2}$ and the norm of $u_{\bel}'$ tends to
zero in $\C^{2,\alpha}_{\nu}$, but the norm of $v_{\bel}'$
in $\C^{4,\alpha}_{\nu+2}$ stays bounded.  If the supremum
of $r^{-\nu-2} |v_{\bel}'|$ occurs at some point $x_{\ell}$, then
exactly the same arguments as above show that $|x_{\ell}|$ stays bounded
away from $0$ and $\infty$. Hence once again we could pass to a limit
to get a solution $V \in \cal C^{4,\alpha}_{\nu+2}(\RNO)$
such that $L_1^* V = 0$. But we know that $L_1^*$ has no elements
in its nullspace which decay at infinity, which implies that $V = 0$,
a contradiction.

Finally, then, we have arrived at the situation above, that
there exist $u \in \C^{2,\alpha}_{\nu}(\RNO)$
and $v \in \C^{4,\alpha}_{\nu+2}(\RNO)$ such that $L_1 u = 0$, $L_1^* v = u$.
Combining these two equations, we get
$L_1 L_1^* v = 0$. But since $v$ decays at infinity, we can
multiply this equation by $v$ and integrate. Integration
by parts now shows that $L_1^* v = 0$, which we have already
observed implies that $v = 0$, hence $u = 0$ as well. This
is the final contradiction. The only alternative is that
the $\C^{2,\alpha}_{\nu}$ norm of $\ubel$ can not blow up.

The case when some $\Sigma_i$ are of positive dimension is treated very
similarly. In fact, only the rescaling needs to be done slightly
differently, but in the same manner as at the end of the proof of
Theorem \ref{th:inj-omega}. \hfill $\Box$

\section{The fixed point argument}

We are now in a position to complete most of the proofs of
Theorem~\ref{th:ex-is-sing} and Theorem~\ref{th:ex-gen-sing}; what
will remain after this section is the assertions about the moduli
of solutions. We need to find a solution $w \in \C^{2,\alpha}_{\nu,\cal D}
(\OS)$ to the equation
\begin{equation}
\Lbe v + Q(v) +\fbe =0
\label{eq:perturbation}
\end{equation}
where $f_{\be}$ is the error term introduced in \S4, either for the
case when $\Sigma$ is discrete, or in the more general case, and $Q$ is the
quadratic remainder (\ref{eq:quad-remain}).
We do this by the standard contraction mapping argument. We will define
a continuous operator $\cal K$ from the space $\C_{\nu,\cal D}^{2,\alpha }
(\OS)$ into itself and then prove that this map
is a contraction on some small ball in this space.

If $v$ is a solution of (\ref{eq:perturbation}), then $\bue +v$
is a weak solution of
\begin{equation}
\left\{ \begin{array}{ll}
-\Delta(\bue+v) = |\bue+v|^p
\quad \mbox{in}\quad  \Omega \\[5mm]
   \bue+v =0 \quad \mbox{on} \quad \del  \Omega.
\end{array}
\right.
\end{equation}
Assuming a solution exists, let us show that $\bue+v$ is positive in
$\OS$. On the one hand, for $x$ near $\Sigma_i$, there exists some $R > 0$
such that if $\rho(x) \le R\e_i$ then
\begin{equation}
c_1 \rho(x)^{-\frac{2}{p-1}} \le \bue(x) \le c_2 \rho(x)^{-\frac{2}{p-1}}
\label{eq:u-not}
\end{equation}
On the other hand, since $v\in \C_{\nu}^{2,\alpha}(\OS)$, we have
$v(x) \leq c \rho(x)^{\nu}$. But since $\nu > -2/(p-1)$, it follows that
$ \bue+v >0$ near each $\Sigma_i$; by the maximum principle we see that
$\bue+v>0$ in all of $\OS$, and hence $\bue+w$ is a positive solution of
(\ref{eq:delta-u=u-p}) which is singular at all points of $\Sigma$.

It remains to prove the existence of a solution to (\ref{eq:perturbation}).
We shall first treat the case where $\Sigma$ is finite.
First observe that $||\fbe||_{0,\alpha,\nu - 2} \le C \e_0^{N -
\frac{2p}{p-1}}$ by Lemma \ref{le:estimate-1}. Let $A$ denote a common
upper bound for the norm of $\Gbe$ for $\e_0$ sufficiently small. Then
\begin{equation}
|| \Gbe \fbe ||_{2,\alpha,\nu} \le AC \e_0^{N - \frac{2p}{p-1}}.
\label{eq:est-remainterm}
\end{equation}

In view of this estimate, we shall work in the ball
\[
{\cal B}(\e_0,\beta) \equiv \{v \in \C_{\nu}^{2,\alpha}(\OS): ||v||_{2,
\alpha,\nu} \leq \beta \e_0^{N - {{2p}\over {p-1}}} \}.
\]
We have already shown that $\Gbe (\fbe)\in {\cal B}(\e_0, \beta)$, for
$\beta$ large enough. At this point we shall fix $\beta$ large enough so
that $\Gbe \fbe \in \cal B(\e_0,\beta/2)$. We shall now make a further
restriction on the $\e_i$, namely that, for some fixed constant $a \in
(0,1)$, we have
\begin{equation}
a \e_0 \le  \e_{i} \le \e_0
\label{eq:cone-isolated}
\end{equation}
for $i = 1, \dots, K$. This new restriction appears
to be needed in the proof of the following Lemma.
\begin{lemma}
There exists some constant $c>0$ independent of $\e_0 \ll 1$ such that
\[
|| Q(v_2) - Q(v_1)||_{0,\alpha,\nu-2} \leq  \frac{1}{2A}\,
||v_{2}-v_{1}||_{2,\alpha,\nu},
\]
for all $v_{1}, v_{2} \in {\cal B}(\e_0,\beta)$. In particular, taking
$v_2 = 0$ we see that $\Gbe Q(v) \in \cal B(\e_0,\beta/2)$ and hence
that the operator $\cal K$ defined by $\cal K(v) = -\Gbe(Q(v) + \fbe)$
maps $\cal B(\e_0,\beta)$ to itself, and is a contraction on this ball
for $\e_0$ sufficiently small.
\label{le:contraction}
\end{lemma}
{\bf Proof~:}
We first establish that there exists some $\tau > 0$, independent of
$\e_0 \ll 1$, such that for any $v\in {\cal B}(\e_0,\beta)$ we have
\[
x \in \cup_{i=1}^K B(x_i,\tau) \Longrightarrow |v(x)| \leq {1\over 4}\bue(x).
\]
Indeed, by Lemma~\ref{pr:base}, we know that there exist constants
$c_1, c_2 >0$ and a radius $R>0$ such that
\[
c_1 |x|^{-\frac{2}{p-1}} \le u_{\e_i}(x) \le c_2 |x|^{-\frac{2}{p-1}}
\quad \mbox{if} \quad |x|\leq R\e_i,
\]
\[
c_1 \e_i^{N - \frac{2p}{p-1}} |x|^{2-N} \le u_{\e_i}(x) \le  c_2 \e_i^{N-
\frac{2p}{p-1}}|x|^{2-N} \quad \mbox{if} \quad R\e_i \leq |x| \le \tau.
\]
The claim follows at once from these estimates and the fact that
$v \in \cal B(\e_0,\beta)$, so that
\[
|v(x)| \le \beta \e_0^{N - \frac{2p}{p-1}} \rho^{\nu}(x) \le c \beta \e_i^{N -
\frac{2p}{p-1}} \rho^{\nu}(x)
\]
by (\ref{eq:cone-isolated}).  We also note that
\[
\bue(x) \le c \rho(x)^{-\frac{2}{p-1}}
\]
for all $x \in \Omega$.

Since $|v/\bue| \le 1/4$ in each $B(x_i,\tau)$, we may use
a Taylor expansion to obtain
\[
| Q(v_2) - Q(v_1)|(x) \leq  c |\bue|^{p-2} (x)|v_2 -v_1|^2 (x).
\]
So for $x\in B(x_{i}, \tau)$ we have
\[
\rho(x)^{2 - \nu} | Q(v_2) - Q(v_1)|(x) \leq c \rho(x)^{ 2 -\nu
-2{{p-2}\over {p-1}} + 2\nu} ||v_{2}-v_{1}||^2_{ 2, \alpha ,\nu}
\]
\[
= c \rho(x)^{\nu +{2 \over {p-1}}} ||v_{2}-v_{1}||^2_{2,\alpha,\nu} \leq c
\tau^{\nu +{2 \over {p-1}}} \beta \e_0^{N - \frac{2}{p-1}}
||v_{2}-v_{1}||_{2,\alpha,\nu}.
\]
The coefficient here may be taken as small as desired by choosing
$\e_0$ sufficiently small. Outside the union of these balls we use the
estimates
\[
\bue(x) \leq c \e_0^{N-{{2p}\over {p-1}}}
\quad \mbox{and also }\quad
|v(x) |\leq c \e_0^{N-{{2p}\over {p-1}}},
\]
where the constant $c>0$ depends on $\tau$ but not on $v$. For $\rho(x) \ge
\tau$ we can neglect all factors involving $\rho(x)$, hence for all $x \in
\Omega_{\tau}$ we have
\[
\rho(x)^{2 -\nu} | Q(v_2) - Q(v_1)|(x) \le C (\bue^{p-1} + |v|^{p-1})|v_2 -v_1|
\]
\[
\le C \e_0^{(p-1)(N - \frac{2p}{p-1})}|v_2 -v_1| \le \frac{\beta}{BA}
\e_0^{(N-2)p -N}||v_2 -v_1||_{2, \alpha, \nu}
\]
for any constant $B > 0$, provided $\e_0$ is chosen small enough.

We also need to estimate the H\"older norm of $Q(v_1) - Q(v_2)$.
By the fourth part of Lemma \ref{le:fcnspacefacts}, it suffices
to estimate $\rho^{3-\nu}|\nabla(Q(v_1) - Q(v_2))|$. We shall only
sketch this briefly. First, in each $B(x_i,\tau)$ we compute
that
\[
\nabla Q(v) = p \left( (\bue + v)^{p-1} - \bue^{p-1} - (p-1) \bue^{p-1}v
\right) (\nabla \bue) + p\left((\bue + v)^{p-1} -
\bue^{p-1})(\nabla v)\right).
\]
For $x \in \Omega_{\tau}$ there is a similar expression, except that
we must be more careful about the absolute value signs. However, in
$\Omega_{\tau}$ we estimate each term individually and conclude that
\[
\rho^{3-\nu}|\nabla Q(v_1) - \nabla Q(v_2)| \le
C \e_0^{(N-2)p -N-1}||v_2 -v_1||_{2, \alpha, \nu}
\]
there. As before, this provides the proper estimate in this set.
In $B(x_i,\tau)$ we use a Taylor expansion in the expression above
to conclude the proper estimate. This completes the proof of
the Lemma when $\Sigma$ is discrete.

To carry out the proof in the more general case, only fairly minor
changes need to be made in the argument above. The most important
one, already observed in \S 4, is that the weight parameter $\nu$
must now be constrained to lie in the smaller interval
\begin{equation}
\frac{-2}{p-1} < \nu < \min\{\frac{-2}{p-1} + 1 = \frac{p-3}{p-1}, \Re(
\gamma_0^- ) \}.
\label{eq:nu-weight-gen}
\end{equation}
Furthermore, the estimate in Lemma \ref{le:Estimate-1} is slightly
weaker than that in Lemma \ref{le:estimate-1}, hence we must replace
the exponent $N - \frac{2p}{p-1}$ in the above argument by
\[
q = \frac{p-3}{p-1} - \nu.
\]
Thus we show that $\cal K$ is a contraction of the balls
\[
\cal B(\e_0,\beta) = \{v \in \C^{2,\alpha}_{\nu,\cal D}(\OS):
||v||_{2,\alpha,\nu} \le \beta \e_0^q\}.
\]
The details of the argument need very little change since the
nonlinear term $Q(v)$ is the same as before. This completes
the proof in all cases. \hfill $\Box$
\smallskip

Using this Lemma and the remarks preceding it, we now see that
there exists a unique solution $v$ of the equation
\[
v=-{\Gbe} \left( (|v+\bue|^p -\bue^{p}-p\bue^{p-1}w)+\fbe \right)
\]
in $\cal B(\e_0,\beta)$ for every $\e_0$ small enough, and for
every $K$-tuple $\be$ satisfying (\ref{eq:cone-isolated}). Fixing some
$\e_0$ which is small enough for this argument to work, and replacing
$\e_0$ by $\e \le \e_0$ in this whole argument shows that we can find
a solution which is a small perturbation of $\bue$ for every
$\e$ and $\bue$ for which $a\e \le \e_i \le \e \le \e_0$.

\section{The proofs completed}

We have now completed the proof of the existence of solutions for Theorems
\ref{th:ex-is-sing} and \ref{th:ex-gen-sing}. In this section we indicate
the arguments leading to the deformation space of solutions for either
problem. We also indicate the small changes needed to complete the proof of
Theorem 3.

\subsection{The deformation spaces}

Included in the statements of Theorems \ref{th:ex-is-sing} and
\ref{th:ex-gen-sing} are the assertions that the space of all solutions to
the equation (\ref{eq:delta-u=u-p}) on a domain $\Omega \subset \RR^N$
with isolated singularities at $\{x_1, \dots, x_K\} \subset \Omega$ is $K$
dimensional, and in general that if any component $\Sigma_i$ of the
singular set is positive dimensional, then there is an infinite dimensional
family of solutions with the same singular set. We recover both of these
statements using the implicit function theorem. The latter case is
analogous to the one studied in \cite{M-S}, and the former mirrors the
situation in \cite{M-P-U}.

In both cases, the main point is that if $u$ is the positive solution to
(\ref{eq:delta-u=u-p}) constructed by the procedure of this paper, then the
linearization $L = \Delta + pu^{p-1}$ to the equation at $u$ is nondegenerate
in the sense that it is surjective as a map from $\C^{2,\alpha}_{\nu,\cal D}
(\OS)$ to $\C^{0,\alpha}_{\nu-2}(\OS)$. This means, using the
implicit function theorem, that all solutions to
the equation which can be obtained as perturbations off of $u$ by terms
growing no faster than $\rho^{\nu}$ are parametrized locally by elements
of the nullspace of this linearization in $\C^{2,\alpha}_{\nu,\cal D}$.
We shall prove this surjectivity in the next proposition, or rather,
we will prove the dual statement as in \S 6 that the linearization
is injective on $\C^{2,\alpha}_{\mu,\cal D}$.

\begin{proposition} If $u = \bue +v$ is a solution of (\ref{eq:delta-u=u-p}),
as constructed in the last section, with all $\e_i \in (a\e,\e)$, and
$v \in \C^{2,\alpha}_{\nu,\cal D}(\OS)\cap \cal B(\e,\beta)$, then if
$\e$ is sufficiently small, the linearization $L_u = \Delta + pu^{p-1}$
is injective as a map on $\C^{2,\alpha}_{\mu,\cal D}(\OS)$.
\label{pr:sol-nondeg}
\end{proposition}
\smallskip
{\bf Proof:} Suppose not, i.e. suppose there exists some sequence $\bel$ with
$\el_i \in (a\el,\el)$ for some sequence $\el \rightarrow 0$, and
a solution $u_{(\ell)} = \bubel + v$ such that the corresponding
linearization $L_{u_{(\ell)}}$ is not injective. Thus there exists
some $\phi^{(\ell)} \in \C^{2,\alpha}_{\mu,\cal D}$ such that
$L_{u_{(\ell)}} \phi^{(\ell)} = 0$. We make exactly the same rescalings
as in the proof of Theorem \ref{th:inj-omega}, to arrive at the
same contradiction. (Note, however, that since the $\el_i$ must all
tend to zero simultaneously, the proof simplifies somewhat.) The only point
that needs checking is that the operator obtained in the limit of these
rescalings is just $L_1$ or $\LL_1$, according as whether $\Sigma$ is
discrete or not. The rescaling sending $u_{\e}(r)$ to $u_1(r)$ is
$u_1(r) = \e^{2/(p-1)}u_{\e}(\e r)$. Applying this same rescaling
to $v(r)$ (suitably translated so that the origin corresponds to
the appropriate point of $\Sigma$) yields the function
$v_{\e}(r) = \e^{2/(p-1)}v(\e r)$. Since $|v(r)| \le \beta \e^q \rho^{\nu}$
for some positive exponent $q$, we see that $|v_{\e}(r)| \le
\beta \e^{q + \nu + \frac{2}{p-1}} \rho^{\nu}$ in some ball $B_{\sigma / \e}$.
Since $\nu + \frac{2}{p-1} > 0$ and $q > 0$, this term tends to zero
with $\e$. Hence the limiting operator is just $\Delta + pu_1^{p-1}$,
as desired. As noted above, this is a contradiction, so the proposition
is proved. \hfill $\Box$
\smallskip

The next step is to compute the size of the nullspace of $L_u$ on
$\C^{2,\alpha}_{\nu,\cal D}(\OS)$. We have already indicated that when some
$\Sigma_i$ is positive dimensional this nullspace is infinite dimensional,
and hence there is an infinite dimensional space of solutions to the
nonlinear problem. But when $\Sigma$ is a finite set, then $L$ is
Fredholm and the nullspace is finite dimensional.  We may calculate
its dimension using a relative index theorem.

\begin{proposition} If $\Sigma$ is discrete, this nullspace of $L_u$ is
$K$ dimensional. If any $\Sigma_i$ has positive dimension, then the
nullspace is infinite dimensional.
\end{proposition}
\smallskip
{\bf Proof:} The latter statement, that the nullspace is infinite dimensional
when some $\dim \Sigma_i > 0$ follows from the theory of \cite{M2}. In
fact, this nullspace is parametrized by functions (of a certain negative
distributional order!) along each of positive dimensional components of
$\Sigma$. This infinite dimensionality can also be seen as a consequence
of the fact that the globalized linearization $\LL_1$ from \S 6.4
has infinite dimensional nullspace, as was pointed out in Proposition
\ref{pr:surj-glob}.

That the nullspace of $L_u$ on $\C^{2,\alpha}_{\nu,\cal D}$ is $K$
dimensional when $\Sigma$ is discrete is proved in exactly the same
manner as the analogous fact in \cite{M-P-U}. We sketch this argument
briefly here. Once again, $L^2$ techniques are used. Let $\text{ind\,}
(\delta)$ denote the index of $L_u$ on $\rho^{\delta}L^2$. From the
self-adjointness of $L_u$ on $L^2(\OS)$ (with respect to Euclidean
measure), we find that $\text{ind\,}(-\delta) = -\text{ind\,}(\delta)$,
for any $\delta \notin \{\delta_j^{\pm}\}$. Although the
index itself requires global information to calculate, it is well
known that the `relative index,' i.e. the difference $\text{ind\,}(\delta')
- \text{ind\,}(\delta'')$ depends only on asymptotic data at the points of
$\Sigma$. Using the relative index theorem proved by Melrose, we see that
the relative index $\text{ind\,}(-\delta) - \text{ind\,}(\delta)$ equals
$2K$ (the details are written out carefully in \cite{M-P-U}. Combining these
two facts, we find that $\text{ind\,}(-\delta) = K$. However, since we
have proved that $L_u$ has no cokernel as a map from $\C^{2,\alpha}_{\nu,
\cal D}$, we conclude that the nullspace of $L_u$ is $K$ dimensional.
\hfill $\Box$
\smallskip

Using these facts, and a standard implicit function theorem argument
as in \cite{M-P-U}, we conclude the following result:
\begin{theorem} When $\Sigma$ is a discrete set, then the solution
space to the equation (\ref{eq:delta-u=u-p}) is locally a $K$ dimensional
real analytic variety. The solutions constructed in this paper lie
in the smooth set of this variety.
\end{theorem}
The fact that this solution space is a real analytic variety (possibly
with singularities) may be deduced using the Ljapunov-Schmidt reduction
argument, as in \cite{K-M-P}. But since we have produced smooth
points in this variety, having found solutions $u$ for which the
corresponding linearization $L_u$ is surjective, we conclude that
the top stratum of this variety is $K$ dimensional, hence that almost
every solution is nondegenerate in this sense.

Finally, we may sharpen the deformation result when $\Sigma$ is not
necessarily discrete, applying the implicit function argument from \cite{M-S}
to conclude
\begin{theorem}Let $\Sigma \subset \Omega$ be any union of $\C^{3,\alpha}$
submanifolds $\Sigma_i$ of dimensions $k_i$ satisfying the restrictions of
Theorem \ref{th:ex-gen-sing}. Then the equation (\ref{eq:delta-u=u-p})
has an infinite dimensional family of solutions.
\end{theorem}

\subsection{The singular Yamabe problem on manifolds of positive
scalar curvature}

The modifications in the arguments of this paper required to solve the
singular Yamabe problem on an arbitary compact manifold $(M,g_0)$, where
$R(g_0) \ge 0$, with singular set prescribed on an arbitrary finite disjoint
union of smooth submanifolds $\Sigma_i$ of dimensions greater than zero and
less than (or equal to) $(n-2)/2$ are very minor. The equation that now
must be solved is
\[
\Delta_{g_0} v - \frac{n-2}{4(n-1)} R(g_0) v + \frac{n-2}{4(n-1)} v^{
\frac{n+2}{n-2}} = 0,
\]
\[
v > 0 \ \text{ on } \ M \setminus \Sigma,\qquad \text{sing\,}(v) = \Sigma.
\]
The operator in the first two terms here, i.e. the linear part, is called
the conformal Laplacian $L_{g_0}$ associated to the metric $g_0$.  The
linearization of this operator around the approximate solution $\bue$ is
\[
\Delta - \frac{n-2}{4(n-1)}( R(g_0) - p \bue^{p-1}).
\]

Now, the rescalings of this linearization may be effected in exactly the same
ways in local Fermi coordinate systems around the submanifolds $\Sigma_i$.
The extra term, of order zero, in this operator disappears in these
rescalings; the resulting model operator we need to study is exactly the
operator $\LL_1$ from \S6. Thus we prove, as before, that this
linearization is uniformly surjective, provided all components of $\be$
are small enough.  The fixed point argument shows that (\ref{eq:perturbation})
has a solution $v$. The only point of the whole argument that needs
special comment is to indicate where we use the assumption that the conformal
class of the metric $g_0$ is positive. This is when we show that $\bue + v$
remains positive on all of $M \setminus \Sigma$; at this step we use that
$\Delta - \frac{n-2}{4(n-1)}R(g_0)$ satisfies the maximum principle --
and this is only true if $R(g_0)$ is nonnegative. This completes
the proof of Theorem \ref{th:sing-Yam}.

\pagebreak

\end{document}